\documentclass{article}

\usepackage{amssymb}
\usepackage{latexsym}
\usepackage[mathscr]{euscript}

\usepackage{amsmath}
\usepackage{amssymb}
\usepackage{fullpage}

\usepackage{amssymb}
\usepackage{latexsym}
\usepackage[mathscr]{euscript}
\usepackage{pb-diagram}
\usepackage{latexsym}
\usepackage{epsfig}

\newcommand\commentout[1]{\marginpar{\tiny $\backslash$commentout}}

\newcommand\qed{\hfill$\square$}

\def\column#1#2{\mathrel{\mathop{#1}\limits_{#2}}}
\def\holim#1#2{\mathrel{\column{holim}{\mathord{\column{\longleftarrow}{#1}}}#2}}

\def\smashover#1{\column{\wedge}{#1}}
\def\compcirc {\mbox{\hspace{.05cm}}\raisebox{.04cm}{\tiny  {$\circ$ }}}

\newtheorem{Lemma}{Lemma}[section]
\newtheorem{Theorem}[Lemma]{Theorem}
\newtheorem{Proposition}[Lemma]{Proposition}
\newtheorem{Definition}[Lemma]{Definition}
\newtheorem{Corollary}[Lemma]{Corollary}

\newtheorem{Example}{Example}[section]
\newtheorem{Remark}{Remark}[section]

\newenvironment{Proof}{\par\noindent\textbf{Proof}}
{\qed}

\usepackage{graphics}
\title{Project Description}

\title{Derived Completions in Stable Homotopy theory}
\author{Gunnar Carlsson \footnote{Research supported in part by NSF DMS-0406992}  \\ Department of Mathematics, Stanford University  \\Stanford, California 94305}
\begin{document}
\maketitle 
\section{Introduction}  It has long been recognized that the  development of a theory of ring and module spectra, which bears the same relationship to the category of spectra as the ordinary theory of rings and modules does to the category of abelian groups,  is a very desirable thing.  A number of such  theories exist. The different approaches  \cite{May et all} and \cite{Smith et all} solve the problem in a satisfactory way, and more recently (see \cite{orthogonal}) versions using orthogonal spectra and $\Gamma$-spaces are also available.   The goal of transporting  constructions which are available for ordinary rings and modules to this new category of ring spectra is  also a very worthwhile one.   Some of these constructions have already been made by the authors of \cite{May et all} and \cite{Smith et all}. In this paper, we will use the notion of an $S$-algebra (as in \cite{May et all}) as the spectrum version of a ring.    Given an $S$-algebra $A$ and module spectra $M$ and $N$, one can construct spectra $M \smashover{A} N$ and $Hom _A(M,N)$, analogous to the constructions $M \column{\otimes}{A} N$ and $Hom_A(M,N)$ for rings $A$ and modules $M$ and $N$.   From the point of view of topologists, the most important constructions to transport to the category of spectra are those which are {\em homotopy invariant}, i.e. those for which module or ring homomorphisms which induce weak equivalences on underlying spectra induce equivalences on the constructions.  For this reason, one considers only {\em derived constructions}, i.e. constructions which on the algebraic side would always replace a module by  a projective resolution for it, and would replace a ring by a levelwise free simplicial ring.  In the context of these  categories of modules, this means that one replaces rings and modules  by {\em cofibrant} and/or  {\em  fibrant }  objects in the categories of $S$-algebras and module spectra.  Our goal in this paper is to introduce and study the derived version of the completion construction for modules over a commutative ring. 

\noindent We describe the main results.  Let $f: A \rightarrow B$ be a homomorphism of commutative  $S$-algebras, and let $M$ be an $A$-module  spectrum.  Then we define  a cosimplicial $A$-module spectrum ${\cal T}^{\cdot}_A(M;B)$, and  the {\em derived completion of $M$ along the homomorphism $f$ }, denoted by $M^{\wedge}_B$,  is the the total spectrum of  ${\cal T}^{\cdot}_A(M;B)$.  Here are the properties of this construction we will prove in this paper. 

\begin{itemize}
\item{  The construction $M \rightarrow M^{\wedge}_B$ is functorial for homomorphism of module spectra over $A$.}
\item{ The construction is functorial in $B$ in the sense that if $B \rightarrow C$ is a homomorphism of commutative $A$-algebras, we obtain an induced homomorphsm $M^{\wedge}_B \rightarrow M^{\wedge}_C$}
\item{The construction is functorial in $A$ in the following sense.  Let $A \rightarrow B \rightarrow C$ be a diagram of commutative $S$-algebras, and let $M$ be a $B$-module.  Let $\rho_A(M) $ denote the spectrum regarded as an $A$-module spectrum by restriction of scalars along the ring homomorphism $A \rightarrow B$.  Then there is an induced map $\rho _A (M)^{\wedge}_C \rightarrow M^{\wedge}_C$.}
\item{ For any commutative ring $A$, $A$ may be regarded as a commutative $S$-algebra via the Eilenberg-MacLane spectrum $\Bbb{H}(A)$,  and an $A$-module $M$ can be regarded as a module spectrum $\Bbb{H}(M)$. Thus, for a homomorphism of commutative rings $f :  A \rightarrow B$,  it is possible to construct the derived completion $\Bbb{H}(M)^{\wedge}_{\Bbb{H}(B)}$. For general $A$, $B$  and $M$,  it is possible that although $\Bbb{H}(A)$, $\Bbb{H}(B)$,  and $\Bbb{H}(M)$ have no higher homotopy,  the derived completion will.  However, if $A$ is Noetherian, $M$ is finitely generated, and $f$ is surjective, the derived completion of $\Bbb{H}(M)$ coincides with the Eilenberg-MacLane construction $\Bbb{H}(M^{\wedge}_I)$, where $M^{\wedge}_I$ is the usual algebraic completion construction at the ideal $I = Ker(f)$. }
\item{One can show that in a sense, the completion ``depends only on $\pi _0 (B)$ for $(-1)$-connected $S$-algebras''.  The precise statement, which is our Theorem \ref{main}, is that if we are given a diagram $f: A \rightarrow B \rightarrow C$ of (-1)-connected commutative  $S$-algebras with $\pi _0 B \rightarrow \pi _0 C$ an isomorphism, and an $A$-module spectrum $M$, then the map $M^{\wedge}_B  \rightarrow M^{\wedge}_C$ is a weak equivalence of spectra.  }
\item{ For any homomorphism $f: A \rightarrow B$ of (-1)-connected commutative $S$-algebras, and any connective $A$-module spectrum $M$ (i.e. $\pi _sM = 0$ for $s$ sufficiently small), there is a spectral sequence whose $E_2$-term depends only on the structure of $\pi _*M$ as a module over $\pi _0 A$ and of $\pi _0 B$ as an  algebra over $\pi _0 A$, and which converges to $\pi _* M^{\wedge}_B$.  This is our Theorem \ref{algebraictogeometric}}
\end{itemize}

\noindent  By far the most important use of the notion of completion in homotopy has been in $p$-adically or profinitely completing spaces.  These constructions are in a sense ``tame'', in that for spectra with finitely generated homotopy groups, the homotopy groups of the completion can be obtained by algebraically completing the homotopy groups at primes or profinitely, and in general are viewed as a simplification of the homotopy type. When the homotopy groups are not finitely generated, one has situations where there is a single derived functor of completion which contributes to the homotopy groups of the completion.    However, when applied to more complicated   rings, our construction can construct interesting homotopy types from discrete rings, even for finitely generated modules over the ring.  For example, let $A$ be the group ring of the discrete group $\Bbb{Z} /  p^{\infty} \Bbb{Z}  = \bigcup _n \Bbb{Z}/ p^n \Bbb{Z}$, and $\Bbb{F}_p$ as an $A$-algebra via augmentation followed by reduction mod $p$.  If we form the derived completion of $A$ itself as an $A$-module (using the Eilenberg-MacLane construction as above), we obtain the $p$-adically completed group ring on the singular chains on the circle group, regarded as a simplicial ring.  What has in effect happened is that the derived completion construction on a discrete ring has created an $S$-algebra (which is much like  a topological ring),  which coincides with our geometric notion of ``filling in the gaps'' in the group $\Bbb{Z} /  p^{\infty} \Bbb{Z} $, viewed as embedded in the circle group. So in this case, the completion produces interesting homotopy types related to embeddings of the ring within topological rings. This phenomenon has some similarity with the behavior of Quillen's plus construction, which replaces the classifying space of a discrete group with a homologically equivalent space.  In \cite{lawson}, T. Lawson has studied the completion process in the context of the pro-$p$ completion of a finitely generated nilpotent group $\Gamma$, and has shown that the homotopy groups of the completion are strongly related with the homology groups of stable (with respect to dimension) representation varieties for $\Gamma$.        

\noindent We are developing this material for use in applications in algebraic $K$-theory, specifically in order to understand the descent problem for the $K$-theory of fields.  The goal is to obtain a homotopy theoretic model for the $K$-theory spectrum of a field $F$ which depends only on the absolute Galois group $G_F$ of the $F$.  It turns out that a model which is often correct can be constructed out of the $S$-algebra associated  to the symmetric monoidal category of finite dimensional complex representations of $G_F$ precisely by performing the derived completion described in this paper.  Although this is our main application, we hope and expect that our constructions will find use in other contexts. 

\noindent At least three  other notions of derived completion have appeared in the literature, in \cite{greenlees}, \cite{DGI}, and \cite{bousfieldpaper}.  The goal of this paper is to record a precise version of the construction which we can use in future work, but we do make some comments on the relationship of the present construction with those in \cite{greenlees} and \cite{DGI}.   

\noindent  The outline of the paper is as follows.  Section 2 develops the preliminary technical material we will require, including material on the theory of $S$-algebras  and module spectra, as well as the theory of Barr and Beck on simplicial and cosimplicial approximations of spectra using ``triple'' or ``monads''.   Section 3 then defines our derived completion, and develops its elementary properties.  In Section 4, we study the behavior of this construction on discrete rings, and show that it coincides with ordinary completion for finitely generated modules over a Noetherian ring. In Section 5, we compare our construction with two other completions, those constructed in \cite{greenlees} and \cite{DGI}. 
Section 6 proves our invariance result for $S$-algebra homomorphisms inducing an isomorphism on $\pi _0$, and Section 7 constructs the spectral sequence discussed above.  Finally, in Section 8, we develop some interesting examples.

\noindent  The author has had valuable conversations with a number of other mathematicians concerning this material.  In particular, discussions  with Bj{\o}rn Dundas, Bill Dwyer, John Greenlees, Mark Hovey, Rick Jardine,  Mike Mandell, J. Peter May, Haynes Miller, and Brooke Shipley have been especially helpful.  
 \section{Preliminaries}
\subsection{$S$-algebras  and module spectra} As mentioned in the introduction, there are a number of constructions of  categories of spectra which  admit a coherently associative and commutative smash product.     One consequence of these constructions  is that one can develop a theory of ``ring spectra" as the category of monoid objects in a category of spectra relative to the smash product. As mentioned in the introduction, we elect to work with $S$-algebras as our notion of ring spectrum. In a similar way, one can define module spectra $M$  over an $S$-algebra  $R$ as spectra equipped with a map $R \wedge M \rightarrow M$ so that the standard algebraic  diagrams commute.   One can also define the notion of a commutative $S$-algebra,  in terms of the commutativity of an obvious diagram.  
The commutative $S$-algebras form a category in their own right, we we will denote by $\mbox{Alg}_S$.  More generally, if $A$ is a commutative $S$-algebra, we can also construct the category $\mbox{Alg}_A$ of commutative $A$-algebra spectra. Given any commutative $S$-algebra $A$, one can also  define a category $\mbox{Mod}_A$ of module spectra over $A$.  Note that because $A$ is commutative, we do not have to specify whether the module is a right or left module.  For any commutative $S$-algebra  there exist relative notions of smash products (analogous to tensor products over a ring) and $Hom$-spectra (analogous to Hom-modules in algebra).  See \cite{May et all} or \cite{Smith et all} for the particulars of these theories.   We will henceforth work with the $S$-algebra version of this theory as constructed in \cite{May et all}.  The following proposition summarizes the  properties of the categories $\mbox{Alg}_A$ and $\mbox{Mod}_A$ which we will need.  The results can all be found in \cite{May et all}, pp. 140-148.

\begin{Proposition} \label{salgsummary} For any commutative $S$-algebra $A$, the categories $\mbox{\em Alg}_A$ and $\mbox{\em Mod} _A$ can both be equipped with the structure of a Quillen model category (see \cite{dwyerone}) with the following properties. 
\begin{enumerate}
\item{ A morphism in $\mbox{\em Alg}_A$ or $\mbox{\em Mod}_A$ is a weak equivalence if and only if its underlying map of spectra is one.  }
\item{All objects in $\mbox{\em Alg}_A$ or $\mbox{\em Mod}_A$ are fibrant. }
\item{ In each category, there is a functorial way to replace each morphism $f$ with a decomposition $f = p \compcirc i$ with $p$ is a fibration and $i$ a cofibration which is also a weak equivalence.  Similarly for $p$ a fibration and a weak equivalence and $i$ a cofibration.  In particular, there exist functorial cofibrant and fibrant replacements in $\mbox{\em Mod}_A$ and $\mbox{\em Alg}_A$.  }
\item{ For any homomorphism $A \rightarrow B$ of commutative $S$-algebras, the functor $ M \rightarrow  B \smashover{A}M$ is a triple in the sense of \cite{barrbeck}. See Section \ref{BarrBecksection} below.  }
\item{There exists a triple $\mbox{\em Sym}_A$ on the category $\mbox{\em Mod}_A$ so that a commutative $A$-algebra spectrum is precisely the same thing as an algebra over the triple $\mbox{\em Sym}_A$, in the sense of \cite{barrbeck}.}
\end{enumerate} 
In particular, the category of spectra is equivalent to the category $\mbox{\em Mod}_S$, where $S$ is the sphere spectrum, and so the results apply to the category of spectra. 
\end{Proposition}  

\noindent     The smash products and Hom constructions are homotopy invariant when the argument modules are cofibrant and or fibrant, in the sense made precise in the following proposition. 

\begin{Proposition} \label{hoequivariance} Let $A$ be a commutative $S$-algebra.  Let $f: M \rightarrow M^{\prime}$ be a weak equivalence of cofibrant $A$-module spectra, and let $N$ be any fibrant $A$-module spectrum.  
Then the induced maps 
$$ f \smashover{A} id_N: M \smashover{A}N \rightarrow M^{\prime} \smashover{A} N
$$
and 
$$ Hom_A(f, N) : Hom_A(M^{\prime}, N) \rightarrow Hom _A (M, N)
$$
are both equivalences.  Also, if $M$ is cofibrant, and $g : N \rightarrow N^{\prime}$ is a weak equivalence, then the natural map 
$$ Hom_A(M,g) : Hom_A (M,N) \rightarrow Hom _A(M, N^{\prime})
$$
is a weak equivalence. 
\end{Proposition}

\noindent The smash product and $Hom$ constructions also behave well on cofibrations of module spectra, in the following sense.  

\begin{Proposition} \label{cofibrations}  Let $A$ be a commutative $S$-algebra, and let $f: M \rightarrow M^{\prime}$ be a cofibration  of $A$-module spectra, and let $N$ be an $A$-module spectrum.  The the induced maps $ f \smashover{A} Id_N$ and $Hom_A(f, N)$ are cofibrations and fibrations, respectively. 
\end{Proposition} 

\noindent These results allow one to develop spectral sequences as computational tools. 

\begin{Proposition} \label{Kunneth} Suppose $A$ is a commutative $S$-algebra, $M$ is a cofibrant $A$-module spectrum, and $N$ is an $A$-module spectrum.  Then there is a spectral sequence with $E_2^{pq} = Tor^{\pi _*(A)} _{pq} (\pi _*(M), \pi _*(N))$, converging to $\pi _{p+q}(M \smashover{A}N)$.  The superscript $p$ refers to homological degree, and $q$ refers to internal degree. 
\end{Proposition}

\begin{Corollary} \label{oneconnectivity}  Suppose $A$ is a (-1)-connected commutative $S$-algebra, $M$ is an $s$-connected cofibrant $A$-module spectrum, and $N$ is a $t$-connected $A$-module spectrum.  Then $M \smashover{A}N$ is $(s +t + 1)$-connected. 
\end{Corollary}

 \subsection{Cosimplicial $S$-algebras  and module spectra}
 \noindent Let $A$ be a commutative  $S$-algebras, and let $\mbox{Mod}_A$ denote the category of module spectra over $A$.  Since $\mbox{Mod}_A$ is a Quillen model category the usual notions of homotopy colimits, homotopy limits, the total $A$-module spectrum of a cosimplicial object in $\mbox{Mod}_A$ as well as its finite stage approximations $Tot$, and fibrancy of a cosimplicial $A$-module spectrum all make sense, and they share the properties of the corresponding notions in the category of simplicial sets and the category of spectra. Every cosimplicial object in $\mbox{Mod}_A$ is functorially equivalent to a fibrant one, and we write $(-)_{fib}$ for this functor.  We recall some of the important properties.  

 \begin{Proposition} \label{holim} Let $F : \Delta \rightarrow {\mbox{Mod}}_A$ denote a fibrant cosimplicial object in $\mbox{Mod}_A$.  Then there is a canonical natural equivalence 
 $$ Tot(F) \rightarrow \holim{\Delta}{F}
 $$
 Moreover, if we let $\Delta ^{(n)}$ denote the full subcategory on the objects of cardinality less than or equal to $n+1$, then there is a natural equivalence 
 $$Tot^n(F) \rightarrow \holim{\Delta ^{(n)}}{F}
 $$
 \end{Proposition}
 
 \begin{Proposition}\label{fubinithm}
 Let $F: \underline{C} \times \underline{D} \rightarrow \mbox{Mod}_A$ be a functor, where $\underline{C}$ and $\underline{D}$ are small categories.  Then for any object $c \in \underline{C}$, we may construct the homotopy inverse limit object $F_c = \holim{c \times \underline{D}}{F \mid c \times \underline{D}}$.  The construction $c \rightarrow F_c$ is functorial in $c$, and we obtain a natural equivalence
 $$ \holim{\underline{C} \times \underline{D}}{F} \hspace{.1cm}\cong \hspace{.1cm}\holim{\underline{C}}{F_c}
 $$
 
 \end{Proposition} 
 \begin{Proposition} Suppose we have a sequence $X^{\cdot} \rightarrow Y^{\cdot} \rightarrow Z^{\cdot}$ of fibrant cosimplicial objects in $\mbox{Mod}_A$,  which is levelwise a cofibration sequence.  Then the sequences 
 $$ Tot^n X^{\cdot} \rightarrow Tot^n Y^{\cdot} \rightarrow Tot^nZ^{\cdot}
 $$ 
 and
 $$ Tot X^{\cdot} \rightarrow Tot Y^{\cdot} \rightarrow Tot Z^{\cdot}
$$
 are cofibration sequences up to homotopy in $\mbox{Mod}_A$.
 \end{Proposition}
 
 \begin{Proposition}\label{smashproperty} Let $X^{\cdot}$ be a fibrant  cosimplicial object in $\mbox{\em Mod}_A$, and let $M$ denote an object in $\mbox{\em Mod}_A$.  Then we may form the new cosimplicial object $M \smashover{A} X^{\cdot}$, and there is a natural equivalence in $\mbox{\em Mod}_A$
 $$ M \smashover{A} Tot^n X^{\cdot} \hspace{.1cm} \cong \hspace{.1cm} Tot^n ( (M \smashover{A} X^{\cdot})_{fib})
 $$
 \end{Proposition}
 \begin{Proof}  Follows easily from (see \cite{May et all} or \cite{Smith et all})  the fact that forming smash products over $A$ with a fixed module commutes with homotopy pullbacks of diagrams of the form $X \rightarrow Y \leftarrow Z$.  We omit the details.  
 \end{Proof}
 
 \noindent It will also be useful to use the idea of an {\em augmented cosimplicial object} in a category.  Let $\Delta_{aug} $ denote the category obtained by adjoining to $\Delta$ the single object $\underline{-1}$, which represents the empty set, regarded as a totally ordered set.  This description also makes it clear how to define the morphisms in $\Delta _{aug}$.  By an augmented cosimplicial object in a category $\underline{C}$, we mean a covariant functor from $\Delta _{aug}$ with values in $\underline{C}$. If we have any abelian group valued functor $A$ on $\underline{C}$, we obtain from any cosimplicial object in $\underline{C}$ an augmented cosimplicial abelian group.  As with cosimplicial abelian groups, we may then associate a  cochain complex (starting in codimension -1) to this cosimplicial abelian group. We summarize the  properties of augmented cosimplicial spectra. 
 
  \begin{Proposition}\label{augmented} Suppose we are given an augmented cosimplicial object $X^{\cdot}$ in a category $\underline{C}$, and let $\rho X^{\cdot}$ denote the cosimplicial object obtained by restriction to $\Delta \subseteq \Delta _{aug}$.  Suppose also that $\underline{C}$ is complete, hence admits a notion of ``total object" $Tot \rho X^{\cdot}$.  Then the coface map $\delta ^0 : X^{-1} \rightarrow X^0$ induces a morphism $\theta : X^{-1} \rightarrow Tot(\rho X^{\cdot})$.  
 \end{Proposition}
 The following result gives a criterion which guarantees that the map $\eta$ is a weak equivalence, when the underlying category $\underline{C}$ is the category $\mbox{Mod}_A$. 
 \begin{Proposition} \label{augmentedone}Suppose $X^{\cdot}$ is an augmented cosimplicial object in the model category $\mbox{Mod}_A$ and that the cosimplicial object $\rho X^{\cdot}$ is a fibrant cosimplicial object in $\mbox{Mod}_A$. For each $t$, we obtain an augmented cosimplicial abelian group $\pi _t(X^{\cdot})$, to which we associate a cochain complex $C^*(t)$ as above.  If each of the cochain complexes $C^*(t)$ has  trivial cohomology, then the map $\theta : X^{-1} \rightarrow Tot(\rho X^{\cdot})$ is a weak equivalence in $\mbox{Mod}_A$.  
 \end{Proposition}
 \begin{Proof} Straightforward verification using the homotopy spectral sequence of a cosimplicial space, see \cite{Bousfield}, Chapter X. 
 \end{Proof}

 \noindent Finally, we will require a  comparison theorem for bicosimplicial $A$-modules.  
  
 \begin{Proposition}\label{levelwise}  Suppose that we have a map of bicosimplicial $A$-modules $f^{\cdot \cdot} : X^{\cdot \cdot} \rightarrow Y^{\cdot \cdot}$.   Suppose that for each $ p \geq 0$, the map $Tot(X^{p \cdot}_{fib}) \rightarrow Tot(Y^{p \cdot}_{fib})$ is a weak equivalence. Then the natural map $Tot (\Delta (X^{\cdot \cdot}) _{fib}) \rightarrow Tot (\Delta (Y^{\cdot \cdot}) _{fib})$ is a weak equivalence, where $\Delta$ denotes the diagonal cosimplicial space.  Similarly if we study the levelwise simplicial objects obtained by holding $q$ fixed. 
 \end{Proposition}
 \begin{Proof} Straightforward, and we omit it. 
 \end{Proof}
   \subsection{The theory of Barr and Beck} \label{BarrBecksection}
Our theory of completion will  make use of the theory of the {\em cosimplicial object of a triple}, a notion discussed by Barr and Beck \cite{barrbeck}.  We will need to use some of the comparison theorems from that paper, so we review that theory here.  For the definition of a triple and algebras over a triple, see \cite{weibel}.  

\noindent Suppose a category $\underline{C}$ is equipped with a triple $T$, and $c \in \underline{C}$.  Then we may define {the cosimplicial resolution of $c$ relative to $T$} to be the cosimplicial object ${\cal T}^{\cdot}(c)$ defined by ${\cal T}^k(c) = T^{k+1}(c)$, and where the cofaces and codegeneracies are defined by 

$$ \delta ^s : T^{k+1}(c) \rightarrow T^{k+2}(c) = T^s (\eta (T^{k+1-s}(c))
$$
and 
$$ \sigma ^t: T^{k+1}(c) \rightarrow T^k(c) = T^t (\mu (T^{k-t-1}(c))
$$
Note that ${\cal T}^{\cdot}(c)$ extends canonically to an augmented cosimplicial object ${\cal T}_{aug}^{\cdot}(c)$ by setting ${\cal T}_{aug}^{-1} = c$, and letting $\delta ^0 : {\cal T}_{aug}^{-1} \rightarrow {\cal T}_{aug}^{0}$ be the canonical inclusion $\eta : c \rightarrow T(c)$.  This means that we have a natural map $\theta : c \rightarrow Tot ({\cal T}^{\cdot}(c))$. 

\noindent Suppose now that the category $\underline{C}$ is $\mbox{Mod}_A$.  For any $A$-module spectrum $M$,  will denote by ${\cal T}_{fib}^{\cdot}(M)$ the functorial fibrant replacement for ${\cal T}^{\cdot}(M)$.  We will need criteria which guarantee that the  map $\theta _{fib}(M) : M \rightarrow Tot ({\cal T}_{fib}^{\cdot}(M))$ is a weak equivalence in $\mbox{Mod}_A$, where $\theta _{fib}(M)$ is the composite $ M \stackrel{\theta}{\rightarrow} Tot({\cal T}^{\cdot}(M)) \rightarrow Tot({\cal T}^{\cdot}_{fib}(M))$.  Barr and Beck now prove the following result. 

\begin{Theorem} \label{BarrBeck}Suppose that $M$ is equipped with a $T$-algebra structure.  Then $\theta_{fib}(M)$ is a weak equivalence. 
\end{Theorem}

\noindent We will derive a useful corollary of this result, using the following lemma. 

\begin{Lemma} \label{insertion}
Let $S$ and $T$ be two triples on the model category $\mbox{Mod}_A$, and suppose we are given a natural transformation $\nu: S \rightarrow T$ of triples.  We may construct the bicosimplicial $A$-module ${\cal C}^{\cdot \cdot} $given by ${\cal C}^{pq} =  {\cal T}^{p}_{fib} ({\cal S}^q_{fib}(M))$.  We may also regard the cosimplicial spectrum ${\cal T}_{fib}^{\cdot}$ as a bicosimplicial spectrum, constant in the $q$-direction, and denote it by ${\cal T}_0^{\cdot \cdot} $.   Then the evident bisimplicial map  $\theta ^T: {\cal T}_0^{\cdot \cdot} \rightarrow {\cal C}^{\cdot \cdot} $ induces an equivalence on total $A$-modules.  
\end{Lemma} 

\begin{Proof}
By applying Proposition \ref{levelwise},  it will suffice to prove that each  of the natural maps 
$$ T^p (M) \rightarrow \column{Tot}{q}(T^p S^q(M)_{fib})
$$
is a weak equivalence of $A$-module spectra.  By \ref{augmentedone}, we need only verify that the cohomology of the cochain complex $A^*$ attached to the augmented cosimplicial group 
$\pi _i (T^p ((S^{\cdot}_{fib}(M)_{aug} ))$ is trivial.  However, for each $q \geq 0$, we have the map 
$h^q: T^p(S^q(M)) \rightarrow T^p(S^{q-1}(M))$ given by the composite
$$
\begin{diagram}
\multiply\dgARROWLENGTH by 3
\node{T^p (S^q(M))}\arrow{e,t}{T^p(\nu(S^{q-1}(M))}  \node{T^p (TS^{q-1}(M))} \arrow{e,t}{T^{p-1}(\mu^T(S^{q-1}(M))}\node{T^{p}(S^{q-1}(M))}
\end{diagram}
$$
As in the proof in \cite{barrbeck} of Theorem \ref{BarrBeck}, the operators $\pi _1 (h^q)$ yield a contracting homotopy for the cochain complex $A^*$, which gives the result. 
\end{Proof}

\noindent We also obtain the following. 
\begin{Proposition} \label{BarrBeckcorollary} Let $S$ and $T$ be two triples on the model category $\mbox{Mod}_A$, and suppose we are given a natural transformations  $\nu: S \rightarrow T$ and $\mu: T \rightarrow S$ of triples.   Then the natural map 

$$ Tot{\cal S}^{\cdot}_{fib}(M) \rightarrow Tot{\cal T}^{\cdot}_{fib}(M)
$$
is an equivalence, where ${\cal S}^{\cdot}$ and ${\cal T}^{\cdot}$ denote the cosimplicial resolutions of $M$ associated to $S$ and $T$, respectively. 
\end{Proposition}
\begin{Proof} This follows directly from  Theorem 5.3 in \cite{radulescu}.  It is also a straightforward consequence of the work in \cite{barrbeck} or \cite{bousfieldone}.  \end{Proof}

\noindent We will also find it useful to discuss the {\em simplicial resolution of a $T$-algebra $X$}, for any triple $T$.  

\begin{Definition} \label{simpresolution} Let $T$ be a triple on a category $\underline{C}$, and let $X$ be any $T$-algebra, with structure morphism $\alpha : TX \rightarrow X$.  The {\em simplicial resolution of $X$ relative to $T$} is the simplicial object $T_{\cdot}(X)$, defined on objects by 
$$  T_k(X) = \underbrace{ T \compcirc T \compcirc \cdots \compcirc T}_{k+1 \mbox{ factors }}(X)
\mbox{ \hspace{.5cm} for $k \geq 0$}$$
and on morphisms by 

$$ 
\left\{ \begin{array}{l}
d_i = T^i\mu (T^{k-i-1}(X)) \mbox{ \hspace{.5cm} for $ 0 \leq i \leq k-1$} \\
\\
d_k = T^k (\alpha) \\
\\
s_i = T^{i+1} \eta ( T^{k-i}(X)) \mbox{ \hspace{.5cm} for $0 \leq i \leq k$}
\end{array}
\right .
$$
The structure map $\alpha $ gives a map  $\varepsilon(X)  : | T_{\cdot}(X) | \rightarrow X$, where $X$ denotes the constant simplicial object with value $X$.   
\end{Definition}

\begin{Proposition} \label{A}Let $T$ denote a triple on the category of spectra.  Then for any $T$-algebra spectrum $X$, the map $\varepsilon (X)$ induces an isomorphism on homotopy groups.  
\end{Proposition}
\begin{Proof} The structure map $\alpha : TX \rightarrow X$ provides an extension of the functor $T_{\cdot}$ to the larger category $\Delta ^{op}_{aug}$, which we will denote by $T_{\cdot}^{aug}$.  We'll refer to a  functor from $\Delta ^{op}_{aug}$ to a category $\underline{C}$ as an {\em augmented simplicial object} in $\underline{C}$.  For any augmented simplicial abelian group $A_{\cdot}$, we may construct the associated chain complex $C_*(A_{\cdot})$, with the alternating sum of face maps as boundary operator.  For any augmented spectrum $X_{\cdot}$, we denote by $\rho(X_{\cdot})$ its restriction to the  subcategory $\Delta ^{op} \subseteq \Delta ^{op}_{aug}$,  so $\rho(X_{\cdot})$ is a simplicial spectrum.  The augmented spectrum gives rise to a map $\nu(X_{\cdot}) : \rho(X_{\cdot}) \rightarrow X_{-1}$. by analogy with the cosimplicial case, it is now easy to verify that given any augmented simplicial spectrum $X_{\cdot}$, if it is the case that for each $i$, the chain complexe 
$C_* \pi _s (X_{\cdot})$ has trivial homology, then the map $\nu(X_{\cdot})$ is a weak equivalence.  
For any triple $T$ on the category of spectra, and every integer $t$, the complex $C_*(\pi _tT^{aug}_{\cdot}(X))$ has trivial homology, since the operators 
$$\pi _t \eta(T^{k+1}(X)) : \pi _t T_{k}^{aug}(X) \rightarrow \pi _t T_{k+1}^{aug}(X)
$$
provide a contracting homotopy for it.  Note that the operator even makes sense for $k = -1$, with $T^0(X) = X$.  This gives the result. 
\end{Proof}

\section{Definitions} \label{definitions}
Let $A$ denote a commutative  $S$-algebra and $B$ a commutative $A$-algebra spectrum. The construction in \cite{May et all} and \cite{Smith et all} of smash products over $A$ make the functor $T_A(-;B) = B \smashover{A}-$ into a triple on the category $\mbox{Mod}_A$. Consequently, we may  construct the cosimplicial resolution of $T_A(-;B)$ as a functor on the category of $A$-module spectra to the category of cosimplicial $A$-module spectra.  We will write ${\cal T}_A^{\cdot}(-;B)$ for the  functorial fibrant replacement of this cosimplicial object in $\mbox{Mod}_A$.

\noindent We recall from \cite{shipley} that for any commutative  $S$-algebra $A$, there is a closed model structure on the category of $A$-algebras.  Moreover, one may functorially  replace any $A$-algebra by a weakly equivalent $A$-algebra which is cofibrant.  

\begin{Definition} Let $A$ denote a commutative $S$-algebra, $B$ a commutative $A$-algebra, and $M$ an  $A$-module.  We define the {\em derived completion} of $M$ at the $A$-algebra $B$ to be the total spectrum of the cosimplicial spectrum ${\cal T}^{\cdot}_A(M;\tilde{B})$, where $\tilde{B}$ is the cofibrant  replacement for $B$ in the category of commutative $A$-algebras.   We write $M^{\wedge}_B$ for $Tot({\cal T}^{\cdot}_A(M;B))$. 
\end{Definition}
The following propositions summarizes the most important properties of this construction.  
\begin{Proposition} \label{properties} Let $A \rightarrow B \rightarrow C$ be a diagram of commutative  $S$-algebras.  Then the following statements all hold.  
\begin{enumerate}
\item{The construction $M \rightarrow M^{\wedge}_B$ is functorial for homomorphisms of  $A$-modules. The map on $B$-completions induced by a homomorphism $f:M \rightarrow N$ of $A$-modules will be denoted by $f^{\wedge}_B$. }
\item{Let $f: M \rightarrow N$ be a homomorphism of  $A$-modules, where $A$ is a commutative $S$-algebra, and let $B$ denote a commutative $A$-algebra.  Suppose that $id_B \column{\wedge}{A} f$ is a weak equivalence of $B$-module spectra.  Then the map on completions $f^{\wedge}_{B}$ is also a weak equivalence of spectra. }
\item{Let $M \rightarrow N \rightarrow P$ be a cofibration sequence of $A$-module spectra,where $A$ is a commutative  $S$-algebra.  Let $B$ be a commutative $A$-algebra spectrum.  Then the sequence 
$M^{\wedge}_B \longrightarrow N_B^{\wedge} \longrightarrow P_B^{\wedge}$ is a cofibration sequence up to homotopy. 
}
\item{Let $A \rightarrow B \rightarrow C$ be a diagram of commutative  $S$-algebras, and let $M$ denote a $B$-module spectrum.  $M$ may be regarded as an  $A$-module spectrum $M^A$.  Suppose that the natural map $C \column{\wedge}{A} M^A \rightarrow C \column{\wedge}{B} M$ is a weak equivalence of spectra.  Then the natural map $(M^A)^{\wedge}_C \rightarrow M^{\wedge}_C$ is also an equivalence of spectra. }
\item{There is a natural transformation $\eta: M \rightarrow M^{\wedge}_B$. }
\item{Suppose $A \rightarrow B$ is a homomorphism of commutative  $S$-algebras.  Suppose $M$ is an $A$-module spectrum, for which the $A$-module structure admits a $B$-module structure extending the given $A$-module structure.  Then the natural map $ \eta: M \rightarrow M_B^{\wedge}$ is an equivalence of spectra. }
\end{enumerate}
\end{Proposition} 
\begin{Proof}
(1) is immediate from the constructions.  (2) is an immediate consequence of the standard fact that a levelwise equivalence of fibrant  spectra induces a weak equivalence of Total spectra.  (3) and  (4) follow from the fact that smash products with a fixed module preserves homotopy cofibration sequences, and that the total spectrum carries levelwise homotopy cofibration sequences to homotopy cofibration sequences.  (5) follows from \ref{BarrBeck}.  (6) is a general property of the cosimplicial resolution of a triple. 
\end{Proof}

\begin{Remark} \label{filtered}{\em We note that although we define the derived completion of an $A$-modules spectrum $M$ at a homomorphism $f: A \rightarrow B$ of commutative  $S$-algebras  to be the total spectrum of ${\cal T}_A^{\cdot}(M;B)$, it is useful to recall that it is therefore the total space of the tower of fibrations 
$$ \cdots \rightarrow Tot^n ({\cal T} ^{\cdot}_A(M;B)) \rightarrow Tot^{n-1} ({\cal T} ^{\cdot}_A(M;B)) \rightarrow \cdots \rightarrow Tot^1 ({\cal T} ^{\cdot}_A(M;B)) \rightarrow Tot^0 ({\cal T} ^{\cdot}_A(M;B))
$$
and that the tower of fibrations is actually functorial as well.  This tower of fibrations gives filtrations on the homotopy groups and other invariants, and will likely be interesting in future $K$-theoretic applications of these ideas. }
\end{Remark}

%For a homomorphism $A \rightarrow B$ of commutative symmetric $S$-algebras,  and an $A$-%module $M$, $A^{\wedge}_B$ is a commutative $A$-algebra, and $M_B$ is a left %$A^{\wedge}_B$ -module. 
%\end{Proposition}

%\begin{Proposition}  Suppose we have a homomorphism of commutative $A$-algebras $B %\rightarrow C$, and that we are given a left $A$-module $M$, so $M^{\wedge}_B$ is an %$A^{\wedge}_B$-module, and $C^{\wedge}_B$ is an $A^{\wedge}_B$-algebra.  Then there is a %canonical (in the homotopy category) equivalence 
%$$  M^{\wedge}_C \cong (M^{\wedge}_B)^{\wedge}_{A^{\wedge}_C}
%$$
%\end{Proposition}

\section{The case of rings}

\noindent For any ring $A$, we  may construct the Eilenberg-MacLane spectrum $\Bbb{H}(A)$.  It is direct from the constructions of \cite{May et all} and \cite{Smith et all} that $\Bbb{H}(A)$ is in a canonical way a  $S$-algebras, and that if the ring is commutative, $\Bbb{H}(A)$ is a commutative  $S$-algebra.  Left and right modules over the  ring $A$  yield, via the $\Bbb{H}$  construction, left and right module spectra over the $S$-algebra $A$.
The goal of this section is to analyze how the derived completion construction given in the preceding section applies to $S$-algebras obtained via this construction.  Specifically, we will show that it coincides with ordinary completion for  finitely generated modules over commutative Noetherian rings, in the sense that $\Bbb{H}(M) ^{\wedge}_{\Bbb{H}(A/I)} \cong \Bbb{H}(M^{\wedge}_I) $ for any finitely generated $A$-module $M$ and any ideal $I \subseteq A$, and where $M^{\wedge}_I$ denotes the usual algebraic completion operation at the ideal $I$.   Note that for a commutative ring $A$, any left $A$-module can be canonically regarded as a right $A$-module as well, so we will simply refer to $A$-modules without specifying right or left modules. 
\begin{Proposition} Let $A$ be a ring, $M$ and $N$   $A$-modules.  Let $\overline{\Bbb{H}(M)} $  denote a cofibrant replacement of the $\Bbb{H}(A)$-module spectrum  $\Bbb{H}(M)$.  Then $$\pi _*(\overline{\Bbb{H}(M)}  \column{\wedge}{\Bbb{H}(A)} \Bbb{H}(N)) \cong Tor^A_*(M,N)$$

\noindent as graded groups.  $M \column{\wedge}{A} N$ denotes the spectrum level construction of ``smash product over $A$'', and Tor denotes the usual algebraic derived functors. 
\end{Proposition}
\begin{Proof}  This results follows  from the K\"{u}nneth spectral sequence whose $E_2$-term is $ Tor^A_*(M,N)$, and which collapses for dimensional reasons.  
\end{Proof}

\noindent  This means that for any family of $k$ left $A$-modules $M_1, M_2, \ldots , M_k$, we may define groups $$\mbox{MultiTor}_i^A(M_1,  M_2, \ldots , M_k)$$ as the  $i$-th homology of the complex $R(M_1) \column{\otimes}{A} R(M_2)  \column{\otimes}{A} \cdots \column{\otimes}{A} 
R(M_k)$, where $R(M_j)$ denotes an $A$-projective resolution  of $M_j$.  

\begin{Proposition} \label{multitor} Let $A$ be a commutative $S$-algebra, and let $M_1, M_2, \ldots , M_k$ denote a family of left $A$-modules.  The K\"{u}nneth spectral sequence generalizes to a spectral sequence with $E_2$-term 
$$ \mbox{\em MultiTor}_*^A(M_1, M_2, \ldots , M_k)
$$
\noindent converging to $\pi _* (M_1 \column{\wedge}{A} M_2 \column{\wedge}{A} \cdots \column{\wedge}{A} M_k )$. 
\end{Proposition}

\begin{Proposition} \label{Bmod} Let $A$ be a commutative ring, and let $B = A/I$, where $I$ is an ideal in $A$.  Suppose further that $M$ is a $B$-module, which  we regard as an $A$-module by restriction of scalars.   Then the natural map $M \rightarrow M^{\wedge}_B$ is an equivalence.  
\end{Proposition}
\begin{Proof} This is immediate from Proposition \ref{properties}, (6). 
\end{Proof}

\noindent What we have now shown is that for any $B$-module $M$, the derived completion $M^{\wedge}_B$ is equivalent to $M$ itself, regarded as a module spectrum over $A$.  We wish to extend this to a result valid for all  finitely generated $A$-modules over a Noetherian ring.  We recall first that the definition of the $I$-adic completion of an $A$-module $M$ is

$$ \column{\mbox{lim}}{\column{\leftarrow}{k}} M/I^kM
$$
\noindent We may also consider the homotopy inverse limit of the pro-$A$-module $\{M/I^kM \}_{k \geq 0}$.  Recall from \cite{Bousfield} that for any inverse system of spectra $\{X_k \}_{k  \geq 0}$, there is a short exact sequence 
$$  0 \rightarrow \mbox{lim}^1 \{  \pi _{i-1}  X_k \} _{k \geq 0} \rightarrow \pi _i \column{\mbox{holim}}{\column{\leftarrow}{k}} X_k \rightarrow \column{\mbox{lim}}{\column{\leftarrow}{k}} \pi _i X_k \rightarrow 0
$$
In our case, the inverse system $\{ \pi _{i-1} X_k \}_{k \geq 0}$ consists entirely of surjective maps, it is a standard result that its ${\mbox{lim}}^1$-term vanishes, leaving us with the isomorphism
$$  \pi _i (\column{holim}{\column{\leftarrow}{k}} M/I^kM) \cong \mbox{\hspace{.01cm}}\column{\mbox{lim}}{\column{\leftarrow}{k}} M / I^kM \mbox{ for $i=0$}
$$
and $ \pi _i (\column{holim}{\column{\leftarrow}{k}} M/I^kM) \cong 0$ for $i \neq 0$.  In other words

$$\column{holim}{\column{\leftarrow}{k}} M/I^kM \cong \Bbb{H}(M^{\wedge}_{I})
$$

\noindent We now wish to show that $M^{\wedge}_B \cong \mbox{\hspace{.01cm}}\column{holim}{\column{\leftarrow}{k}} M/I^kM$.  Consider the inverse system of cosimplicial $A$-module spectra $ \{ {\cal T}^{\cdot}_A(M/I^kM;B) \} _{k \geq 0}$.  We have the natural map 

$$   \theta: {\cal T}^{\cdot}_A(M;B) \rightarrow   \{ {\cal T}^{\cdot}_A(M/I^kM;B) \} _{k \geq 0}
$$
where $ {\cal T}^{\cdot}_A(M;B)$ is considered as a constant pro-cosimplicial $A$-module spectrum.  By taking total spectra and  homotopy inverse limits (in the $k$-direction), we obtain a map of $A$-module spectra $M^{\wedge}_B \rightarrow \column{holim}{\column{\leftarrow}{k}}Tot {\cal T}^{\cdot}_A(M/I^kM ; B)$. By \ref{Bmod}, we have that  $Tot {\cal T}^{\cdot}_A(M/I^kM ; B) \cong M/I^kM$, so we obtain a natural map $\lambda : M^{\wedge}_B \rightarrow \column{holim}{\column{\leftarrow}{k}}M/I^kM \cong \Bbb{H} (M^{\wedge}_I)$.   Our result is now

\begin{Theorem} \label{ringcase} For a Noetherian commutative ring $A$, and a finitely generated $A$-module $M$, the map $\lambda : M^{\wedge}_B \rightarrow  \Bbb{H} (M^{\wedge}_I)$ described above is an equivalence. 
\end{Theorem}
\begin{Proof} It is standard that homotopy inverse limits of cosimplicial spectra commute with taking total spectra, so it is enough to verify that the maps 
$$ {\cal T}^{i}_A(M;B) \rightarrow \column{holim}{\column{\leftarrow}{k}} {\cal T}^{i}_A(M/I^kM ; B)
$$
are equivalences for each $i$.   In order to prove this, it will suffice to show that the pro-$A$-module $\{ {\cal T}^{i}_A(M/I^kM ; B) \}_{k \geq 0}$ is isomorphic, as pro-abelian groups, to the constant pro-abelian group with value ${\cal T} ^{i}_A(M;B)$.  By \ref{multitor}, this means that it will suffice to show that the pro-abelian group 
$$ \{ \mbox{\em MultiTor}_*^A(\underbrace{B, B, \ldots , B}_{s \mbox{ factors}}, M/I^kM) \}_{k \geq 0}
$$
is isomorphic to the constant pro-abelian group with value $ \mbox{\em MultiTor}_*^A(\underbrace{B, B, \ldots , B}_{s \mbox{ factors}}, M)$.    For any $A$-module $M$, we let $\underline{I}M$ denote the pro-$A$-module $\{ I^kM \} _{k \geq 0}$.  Recall (see \cite{adamscompletion}) that the category of pro-$A$-modules is itself an abelian category, and so the notion of exact sequence has meaning.  The construction $M \rightarrow \underline{I}M$ gives a functor from the category of finitely generated $A$-modules to the category of pro-$A$-modules, which is exact in the sense that it carries exact sequences to exact sequences.  The exactness follows from the usual proof of exactness of completion for finitely generated modules over a Noetherian ring (see \cite{Eisenbud}). This same exactness also shows that for complexes of finitely generated $A$-modules, homology commutes with the functor $\underline{I}$, so that we have 
$$\mbox{\em MultiTor}_*^A(\underbrace{B, B, \ldots , B}_{s \mbox{ factors}}, \underline{I}M) \cong \underline{I}\mbox{\em MultiTor}_*^A(\underbrace{B, B, \ldots , B}_{s \mbox{ factors}}, M)
$$ 
It is now clear $\mbox{\em MultiTor}_*^A(\underbrace{B, B, \ldots , B}_{s \mbox{ factors}}, M) $ is a $B$-module, and therefore that $I$ acts trivially on it, which means that the pro-abelian group $\underline{I} \mbox{\em MultiTor}_*^A(\underbrace{B, B, \ldots , B}_{s \mbox{ factors}}, M) $ is pro-isomorphic to the zero group.  We have an exact sequence of pro-abelian groups
$$ \underline{I}M \rightarrow M \rightarrow M / \underline{I}M
$$
so we may conclude that the map $M \rightarrow M/ \underline{I} M$  induces an isomorphism on homology as pro-groups.  This is the required result. 
\end{Proof}

\section{Comparison with other constructions}
In this section we will prove some results comparing the present construction with those of Greenlees-May \cite{greenlees} and Dwyer-Greenlees-Iyengar \cite{DGI}. 
\subsection{The Greenlees-May completion}
This completion construction considers the situation of a commutative  ring $A$ and a finitely generated ideal $I$. It begins with a purely algebraic construction of derived functors of completion, and then produces a spectrum level construction.  We will compare the derived functors defined in \cite{greenlees} with the homotopy groups of our derived completion, and proved that they agree.  It appears likely that the methods of Section \ref{algeo} would allow a comparison for the spectrum level construction, but we do not carry that out here.   We will write $B = A/I$, and let $S$ denote any finite generating set for $I$.  In \cite{greenlees}, the authors construct a chain complex $C_*^S$, essentially a colimit of Koszul complexes based on $S$, and define their derived completion of an $A$-module $M$, which we will denote by $M^{gm}_S$, by 
$$ M^{gm}_S = Hom_A(C^S_*, M)
$$
The construction depends on $S$, but there are canonical isomorphisms between the homology groups constructed using different finite generating sets $S$.  These are the derived functors of completion at $I$ applied to the module $M$.  There is a natural chain map $C^S_* \stackrel{\gamma}{\rightarrow} A_*$, where $A_*$ denotes the chain complex, concentrated in degree 0, with zeroth module equal to $A$.  The map $\gamma$ now induces a natural map 
$$ \xi: M \rightarrow M^{gm}_S
$$
where $M$ denotes the module $M$ regarded as a chain complex concentrated in degree zero. We now have the following. 
\begin{Lemma} \label{BBmod}
Let $M$ be any $B$-module, regarded as an $A$ module by restriction of scalars.  Then the natural map $\xi:M \rightarrow M^{gm}_S$ is a quasi-isomorphism of chain complexes
\end{Lemma}
\begin{Proof} Immediate consequence of Lemma 1.3 of \cite{greenlees}.
\end{Proof}

\noindent It is also immediate from the definitions that the construction $M \rightarrow M^{gm}_S$ extends to chain complexes of $A$-modules.  The following is now easy to prove.
\begin{Lemma} \label{complextest} The construction $C_* \rightarrow (C_*)^{gm}_S$  preserves quasi-isomorphisms. Given a chain complex $C_*$ which is bounded below, and so that   the $A$-module structure on $H_i(C_*)$ extends to a $B$-module structure for all $i$, then the natural map $\xi: C_* \rightarrow (C_*)^{gm}_S$ is a quasi-isomorphism.  
\end{Lemma}
\begin{Proof} The first statement is straightforward from the definitions.  The second one follows from the existence of a ``Postnikov tower" for any bounded below chain complex, where the relative terms are resolutions of the $A$-modules $H_i (C_*)$.  The first statement implies that if $R_*(M)$ is a resolution of an $A$-module $M$, then  the map $R_*(M) \rightarrow M$ induces an equivalence $(R_ *(M))^{gm}_S \rightarrow M^{gm}_S$.  The result now follows from Lemma \ref{BBmod} above, together with the evident fact that the construction carries exact sequences of chain complexes to exact sequences. 
\end{Proof}

\noindent In order to compare $M^{gm}_S$ with $M^{\wedge}_B$, we will need to reinterpret $M^{\wedge}_B$ as a chain complex.  We first begin by constructing a cofibrant version ${\cal B}_*$  of $B$.  This can be done, for example, by using the free commutative $A$-algebra functor on sets.  This is a triple on the category of sets, and for any $A$-algebra yields a simplicial resolution of $B$, which is levelwise free as an $A$-module.  It is also a simplicial commutative ring.  The cosimplicial spectrum ${\cal T}_A(M;{\cal B}_*)$  which constructs $M^{\wedge}_B$ is now a cosimplicial simplicial abelian group, and its total space is a simplicial abelian group.  As such, it corresponds to a chain complex.  On the other hand, one can define the notion of the chain complex associated to  a spectrum. One begins with the associated chain complex functor $C_*$ for the notion of spaces one is dealing with, the singular complex in the case of topological spaces, and the simplicial chains for simplicial sets.  For any spectrum ${\cal X} = \{ X_i \}_{i \geq 0} $, with structure maps $\sigma _i : \Sigma X_i \rightarrow X_{i+1}$, one defines $ch({\cal X})$ to be the colimit of the system 
$$    \cdots \rightarrow \Sigma^{-i}C_*(X_i) \rightarrow \Sigma^{-(i+1)}C_*(X_{i+1}) \rightarrow \Sigma^{-(i+2)}C_*(X_{i+2}) \rightarrow \cdots 
$$
Now, the cosimplicial spectrum corresponds to a cosimplicial chain complex, by applying $ch$  levelwise.  We write $E(k,M)_*$ for the chain complex in codimension $k$.  It is defined by 
$$  E(k,M)_* = ch( \bigwedge _k {\cal B} \wedge M)
$$
For any cosimplicial chain complex $C(-)_*$, we define a chain complex ${\cal D}C(-)_*$ by
$$ {\cal D}C(-)_* = \prod_k \Sigma ^{-k} C(k)_*
$$ 
on the level of graded groups ($\Sigma ^t$ just denotes a shift of  degrees by $t$), and the boundary map comes from the bicomplex structure on this graded group, with one boundary being the one existing on each of the complexes $\Sigma ^{-k}C(k)_*$, and the other being   the alternating sum of the coface maps.  

\begin{Proposition} Let $k \rightarrow A(k)_.$ denote any cosimplicial simplicial abelian group.  Then there are canonical isomorphisms
$$  \pi _s (Tot_k A(k)_.) \cong H_s ( {\cal D} ch(A(-)_.))
$$
for all $s$
\end{Proposition}
\begin{Proof} This is standard.  See for example \cite{Bousfield}. 
\end{Proof}

\noindent We now write ${\cal E}_*(M)$ for the chain complex ${\cal D}E(-,M)_*$.  The construction ${\cal E}_*$ is clearly functorial in $M$, and $H_*({\cal E}_*(M)) $ is canonically isomorphic to 
$ \pi _*(M^{\wedge}_B)$.  We note that ${\cal E}_*$ is a chain complex of $A$ modules, and that the standard map from the constant cosimplicial $A$-module with value $M$ to ${\cal T}_A (M; {\cal B}_*)$ induces a natural transformation $\theta: M_*\rightarrow {\cal E}_*(M)$, where  $M_*$ denotes $M$ regarded as a chain complex concentrated in degree zero.  

\noindent We wish to compare the groups $H_*(M^{gm}_S)$ and $H_*({\cal E}_*(M))$.  The idea will be to compare each of the complexes  in question with the complex ${\cal E}_*(M)^{gm}_S$.  The natural transformation $\theta$ produces a natural map $\Theta: M^{gm}_S \rightarrow {\cal E}_*(M)^{gm}_S$, and the transformation $\xi$ above yields a natural map  $\Xi : {\cal E}_*(M) \rightarrow {\cal E}_*(M)^{gm}_S$.  We now have the following result.  

\begin{Proposition} The natural maps $\Theta$ and $\Xi$ are quasi-isomorphisms.  Therefore the groups $H_*(M^{gm}_S)$ and $H_*({\cal E}_*(M))$ are canonically isomorphic.  
\end{Proposition}
\begin{Proof}  We begin with $\Xi$.  We note first that ${\cal E}_*(M)$ is itself obtained as the total complex of a bicomplex, and can be filtered by the codimension $k$. Note that this is a decreasing filtration.   The subquotients are the complexes $ch(\bigwedge_k {\cal B} \wedge M)$.  We obtain a  filtration on ${\cal E}_*(M)^{gm}_S$  which is compatible with the filtration on ${\cal E}_*(M)$ under the natural map $\Xi$.  Both complexes are obtained as inverse limits of the quotients by the terms in the filtrations, so in order to prove the result it will suffice to show that the natural map 
$$ ch(\bigwedge_k {\cal B} \wedge M) \rightarrow (ch(\bigwedge_k {\cal B} \wedge M))^{gm}_S
$$ 
is a quasi-isomorphism for all $k$.  The reduction uses the $lim^1$ sequences for the homology of inverse limits of chain complexes. But now, in view of Lemma \ref{complextest} above, it will suffice to show that the $A$-module structure on   $H_i (ch(\bigwedge_k {\cal B} \wedge M))$ extends to a $B$-module structure.  But this is clear, since the homology of these complexes are the $MulitTor$ groups of $M$ with coefficients in $B$, and the elements of $B$ act on $B$ by multiplication.  This shows that $\Xi$ is a quasi-isomorphism.  For $\Theta$, we note that the Greenlees-May complex $C_*^S$ is a colimit of subcomplexes $C_*^S(k)$, which are quasi-isomorphic to the Koszul complexes ${\cal K}(k)$  based on the set $\{ s^k \}_{s \in S}$.  It follows that $ M^{gm}_S$ is expressed as the inverse limit of complexes $M^{gm}_S(k)$, which are quasi-isomorphic to $Hom _A({\cal K}(k), M)$, and similarly for chain complexes.  Therefore, there is also  a description of ${\cal E}_* (M)^{gm}_S$ as an inverse limit of quotient complexes ${\cal E}_*(M)^{gm}_S(k)$, which is compatible (under the map $\Theta$) with the inverse limit description of $M^{gm}_S$.  It therefore suffices to show that the natural map 
$$  Hom _A({\cal K}(k), M) \rightarrow Hom _A({\cal K}(k), {\cal E}_*(M))
$$
is a quasi-isomorphism.  To see this, we note that there is a canonical isomorphism of complexes $Hom _A({\cal K}(k), {\cal E}_*(M)) \rightarrow {\cal E}_*(Hom _A ({\cal K}(k), M))$.  We also note that the complex $Hom _A({\cal K}(k), M)$ admits a finite Postnikov tower, with the relative quotients being quasi-isomorphic to resolutions of the various $A$-modules $H_i (Hom _A({\cal K}(k), M))$.  This therefore  induces a similar  tower on  ${\cal E}_*(Hom _A ({\cal K}(k), M))$, and it therefore suffices to show that the natural map $R_*(H_i) \rightarrow {\cal E}_*(R_*(H_i))$ is a quasi-isomorphism for all $i$, where $$H_i = H_i (Hom _A ({\cal K}(k), M))$$ Since a resolution of an $A$-module is quasi-isomorphic to that module regarded as a constant chain complex, we have reduced the problem to verifying that $H_i \rightarrow {\cal E}_*(H_i)$ is a quasi-isomorphism for all $i$, or equivalently that $H_i \rightarrow (H_i)^{\wedge}_B$ is an equivalence.      The action of any element $s^k$ on ${\cal K}(k)$ is easily seen to be chain null homotopic, and the similar result follows for $P = Hom _A({\cal K}(k), M)$.  Each $H_i$ therefore  has the property that every element of the form $s^k$, for $s   \in S$, acts by zero on $H_i$. It follows that there exists an integer $l$ so that the ideal $I^l$ acts trivially on $H_i$. Now consider the filtration $\{I^s H_i \}$ on $H_i$.  It is a finite filtration, and the subquotients $I^s H_i /I^{s+1} H_i$ are $A$-modules so that $I$ acts trivially. Using the exactness property of our derived completion construction, it suffices to prove that the natural maps 
$$ I^s H_i /I^{s+1} H_i \rightarrow (I^s H_i /I^{s+1} H_i)^{\wedge}_B
$$ are equivalences.  But this is Proposition \ref{Bmod}.   \end{Proof}

\noindent The following corollary is now immediate. 
\begin{Corollary} Let $A$ be any commutative ring,  $I$ any finitely generated ideal, with finite generating set $S$.  Let $M$ be any $A$-module.  Then the Greenlees-May completion $M^{gm}_S$ and the derived completion $M^{\wedge}_{A/I}$ are canonically equivalent in the homotopy category of $A$-module spectra, where $A$ is regarded as an $S$-algebra via the Eilenberg-MacLane construction.  
\end{Corollary}
%There is also the cosimplicial chain complex $D(k)_*$ obtained by applying the K\"{u}nneth equivalence to $C(k)_*$, which is defined by 
%$$
%$$
\subsection{The Dwyer-Greenlees-Iyengar  completion} 

\noindent In \cite{DGI}, the authors construct a version of completion which is inspired by Morita theory.  The idea is as follows.  Let $A$ and $B$ be commutative $S$-algebras, and let $f: A \rightarrow B$ be a homomorphism of commutative $S$-algebras.  One can then construct an $S$-algebra ${\cal E} = {\cal E}_A(B) = Hom_A(B,B) = End_A (B)$, where the $Hom$ constructions are in the category of $A$-modules.  $B$ now becomes a left ${\cal E}$ module.  Further, for any $A$-module $M$, we may form $B \column{\wedge}{A}M$, and it becomes a left ${\cal E}$-module  as well.  One now constructs $Hom_{{\cal E}}(B, B \column{\wedge}{A} M)$, and note that there is a natural homomorphism $M \rightarrow Hom_{{\cal E}}(B, B \column{\wedge}{A} M)$, whose adjoint is the map 
$$ B \wedge M \rightarrow B \column{\wedge}{A} M
$$ 
of left ${\cal E}$-modules, where the action on $B \wedge M$ is on the first factor.  We will call the  spectrum $Hom_{{\cal E}}(B, B \column{\wedge}{A} M)  $  the {\em DGI completion} of $A$ at the homomorphism $f$.  The main subject of \cite{DGI} is duality theory in module spectra, and this notion of completion is best understood in the context where $B$ satisfies certain finiteness conditions, which is where the duality theory of \cite{DGI} takes place.  We will see that in our situation, there is a related (more restrictive)  finiteness condition, under which we will prove that the DGI completion agrees with our derived completion.  

\noindent Given $A$,$B$,$M$, $f$, and ${\cal E}$ as above, we construct the cosimplicial $A$-module spectrum ${\cal T}_A^{\cdot}(B;M)$ as in Section \ref{definitions}.  We can now construct the cosimplicial $B$-module spectrum $B \column{\wedge}{A} {\cal T}_A^{\cdot}(B;M)$, to obtain a cosimplicial ${\cal E}$-module ${\cal M}^{\cdot}$, and construct the cosimplicial ${\cal E}$-module spectrum
$$  Hom _{{\cal E}} (B, {\cal M}^{\cdot})
$$ 
Of course, there is the canonical map $M \rightarrow {\cal T}_A^{\cdot}(B;M)$, consequently the induced map $B \column{\wedge}{A} M \rightarrow {\cal M}^{\cdot}$, and therefore a natural map 
$$ \alpha:  Hom _{{\cal E}}(B, B \column{\wedge}{A} M) \rightarrow 
Hom_{{\cal E}}(B,{\cal M}^{\cdot})
$$
where the left hand side is the DGI completion of $M$.  

\begin{Proposition} \label{derdgi}
The map $Tot(\alpha )$ of total spectra induced by $\alpha$ is an equivalence of spectra.  
\end{Proposition}
\begin{Proof} One readily checks that there is an equivalence 
$$ Hom_{{\cal E}}(B,Tot {\cal M}^{\cdot}) \rightarrow Tot(Hom_{{\cal E}}(B,{\cal M}^{\cdot}))
$$
and reduces to proving that the natural map $B \column{\wedge}{A} M \rightarrow Tot({{\cal M}^{\cdot}})$ is an equivalence.  This follows readily from Theorem \ref{BarrBeck}. 
\end{Proof}

\noindent Since we have a natural map from an $A$-module $M$ to its DGI completion, we obtain a natural map 
$$  \rho: {\cal T} _A^{\cdot}(M;B) \rightarrow Hom _{{\cal E}} (B, {\cal M}^{\cdot})
$$
and the corresponding map of total spectra
from $M^{\wedge}_B $ to $Tot (Hom _{{\cal E}} (B, {\cal M}^{\cdot}))$, which in the homotopy category can be interpreted as a map from $M^{\wedge}_B $ to the DGI completion of $M$, according to Proposition \ref{derdgi}.   

\begin{Lemma}  Suppose that $B$ is finitely built from $A$ in the sense of \cite{DGI}, and that the $A$-module $M$ is also finitely built from $A$.  Then the map $\rho$ defined above is an equivalence.  In other words, when $B$ is finitely built from $A$, our derived completion agrees with the DGI completion. 
\end{Lemma} 
\begin{Proof}  We  sketch a proof.  One can write $Hom_{{\cal E}}(B, Tot {\cal M}^{\cdot})$ as the total spectrum of the cosimplicial spectrum 
$$  k \mapsto Hom_{{\cal E}}(B,  B \column{\wedge}{A}{\cal T}^k_A(M;B))
$$
and in order to prove the result it will clearly suffice to prove that the natural maps 
$${\cal T}^k_A(M;B) \rightarrow Hom _{{\cal E}}(B, B \column{\wedge}{A}{\cal T}^k_A(M;B))
$$
are equivalences.  From the definition of ${\cal T}_A^{\cdot}$, it will now suffice to show that for any module $N$ which is finitely built from $A$, the natural map
$$  B \column{\wedge}{A} N \rightarrow 
Hom_{{\cal E}}(B, B \column{\wedge}{A} B \column{\wedge}{A} N)
$$
is an equivalence. This  follows from the three facts  given below, which are easily obtained using the results of \cite{DGI}.  
\begin{itemize}
\item{The correspondence $D$ defined by  $D(N)  =  Hom_A (N,A)$ is a contravariant equivalence of categories from the category of modules which are finitely built over $A$ to itself.  There is a canonical equivalence $N \cong D^2(N)$.  }
\item{Given any $A$-module $N$ which is finitely built from $A$, there is a canonical equivalence 
$$ Hom_{{\cal E}}(B, B \column{\wedge}{A} D(N)) \rightarrow 
Hom_{{\cal E}}(B \column{\wedge}{A} N, B)
$$}
\item{For any $A$-module $N$ which is finitely built from $A$, there is an equivalence of left ${{\cal E}}$-modules 
$$B \column{\wedge}{A} D(B \column{\wedge}{A} N) \cong {\cal E} \column{\wedge}{A} D(N)$$
where the ${\cal E}$-action on the left hand factor is on the left hand factor $B$, and on the right hand side by left multiplication on the ${\cal E}$-factor. }
\end{itemize} 
\end{Proof}

\noindent We now extend this Lemma into the main result of this section.  We recall the notion of {\em proxy finiteness} from \cite{DGI}.  Let $A $ be a   commutative $S$-algebras.  The  $A$-module $B$ is said to be proxy finite if there is a $A$-module $K$ so that $K$ is finitely built from $A$, and so that $K$ is built from $B$.  The DGI completion uses arbitrary $A$-modules as input, while our construction requires a commutative $A$-algebra.  Consequently, this notion as it stands is not useful for us.  We define the following replacement notion which is what is needed in our context. 
\begin{Definition} Let $A$ be a commutative $S$-algebra, and let $B$ be a commutative $A$-algebra, with structure homomorphism $f: A \rightarrow B$.  We say $B$ is {\em algebra proxy finite} over $A$ if there is a commutative diagram 
$$
\begin{diagram}
\node{A} \arrow{e} \arrow{se,b}{f}\node{\overline{B}} \arrow{s}  \\
\node{} \node{B}
\end{diagram}
$$
of commutative $S$-algebras, where the $A$-module  $\overline{B}$ is finitely built from $A$, and built from $B$.  
\end{Definition}

\noindent The following lemma makes this notion useful.  
\begin{Lemma} \label{proxy} Let 
$$
\begin{diagram}
\node{A} \arrow{e}\arrow{se}\node{\overline{B}} \arrow{s}  \\
\node{} \node{B}
\end{diagram}
$$
be a commutative diagram of commutative $S$-algebras as above, and suppose $\overline{B}$ is built from $B$.  Then for any $A$-module $M$, the natural map 
$$  M^{\wedge}_{\overline{B}} \rightarrow M^{\wedge}_B
$$
is an equivalence of $A$-module spectra. 
\end{Lemma}
\begin{Proof}  We consider the bicosimplicial spectrum ${\cal B}^{\cdot \cdot}$ defined by 
$$ (k,l) \mapsto {\cal T}_A^k ( {\cal T}_A^{l}(M;\overline{B});B)
$$
It is equipped with a natural map 
$$\eta:  {\cal T}_A^{\cdot}(M;\overline{B}) \rightarrow {\cal B}^{\cdot \cdot}
$$
when ${\cal T}_A ^{\cdot}(M;\overline{B})$ is regarded as a cosimplicial spectrum constant in the $k$-direction.  We first observe that $\eta$ induces an equivalence on total spectra.  To verify this only requires that we show that for any $\overline{B}$-module $N$, the natural map $N \rightarrow N^{\wedge}_B$ is an equivalence.  This follows from the fact that the total spectrum construction respects cofibration sequences, together with Proposition  \ref{properties}, part 5.  There is also a natural map $\nu: {\cal T}_A^{\cdot}(M;B) \rightarrow {\cal B}^{\cdot \cdot}$, where ${\cal T}_A^{\cdot}$ is regarded as a bicosimplicial spectrum constant in the $l$-direction.  It too induces an equivalence on total spectra.  It is easy to see that to verify this requires only that we show that for any $A$-module $N$, the natural map 
$$ B \column{\wedge}{A} N \rightarrow Tot (B \column{\wedge}{A} {\cal T}_A ^{\cdot}(N;\overline{B}))
$$
is an equivalence of spectra.  But, there is an evident equivalence of cosimplicial spectra
$$B \column{\wedge}{A} {\cal T}^{\cdot}_A (N;\overline{B}) \simeq
{\cal T}^{\cdot}_A(B  \column{\wedge}{A} N; \overline{B})
$$    

\noindent and we are reduced to showing that the map
$$ \eta : B \smashover{A} N \rightarrow Tot {\cal T}^{\cdot}_A(B \smashover{A} N; \overline{B}) \cong (B \smashover{A} N)^{\wedge}_{\overline{B}}
$$
is an equivalence.  But, since $B$ an $\overline{B}$-algebra spectrum, it follows directly that $B \smashover{A} N$ admits a $\overline{B}$-module structure extending the $A$-module structure, and the result now follows directly from 
 Proposition \ref{properties}, part 5.  It is finally not difficult to check that we have a commutative diagram in the homotopy category 
$$  
\begin{diagram}
\node{M^{\wedge}_{\overline{B}}} \arrow{e}  \arrow{se,b}{\eta} \node{M^{\wedge}_B} \arrow{s,t}{\nu} \\
\node{} \node{Tot{\cal B}^{\cdot \cdot}}
\end{diagram}
$$
where the horizontal map is the usual map arising from the functoriality of the completion construction.  This gives the result. 
\end{Proof} 

\noindent We now draw our main conclusion.  
\begin{Theorem}  Let $f: A \rightarrow B$ be a  homomorphism of commutative $S$-algebras, and suppose that $B$ is algebra proxy finite.  Then for any $A$-module $M$ which is finitely built from $A$, the map $\rho $ defined above from $M^{\wedge}_B $ to the DGI completion of $M$ is an equivalence.  
\end{Theorem}
\begin{Proof} Follows directly from Lemma \ref{proxy} above, together with Theorem 4.10 of  \cite{DGI}, from which it follows that the DGI completion has an invariance property similar to that proved in Lemma \ref{proxy} above.  
\end{Proof} 
\section{Main isomorphism theorem}
In this section we will prove the following theorem about derived completions. 

\begin{Theorem}\label{main} Suppose that we have a diagram $A \rightarrow B \stackrel{f}{\rightarrow} C$ of commutative  $S$-algebras. Suppose further that $A,B$, and $C$ are all (-1)-connected, and that the homomorphism $\pi _0 (f)$ is an isomorphism.  Suppose further that the natural homomorphisms $\pi _0(A) \rightarrow \pi _0 (B)$ and $\pi _0 (A) \rightarrow \pi _0 (C)$ are surjections. Then for any left $A$-module spectrum $M$,  the natural homomorphism $M_B \rightarrow M_C$ is an equivalence of $A$-module spectra. 

\end{Theorem}

\noindent The proof of this theorem requires some preliminary technical work on cosimplicial spaces and spectra.   Let $X^{\mbox{\large $\cdot$ }}$ denote any fibrant  cosimplicial space (or, more generally, a fibrant cosimplicial spectrum or $A$-module, where $A$ is a commutative $S$-algebra), i.e. a space-valued functor  from the category $\Delta$ whose objects are the totally ordered sets $\underline{k} = \{ 0,1,\ldots, k \}$ and whose morphisms are the ordering preserving maps of sets.  From Proposition \ref{holim}, we have that  $Tot(X^{\mbox {\large $\cdot$}})$ is weakly equivalent to $\column{holim}{\column{\leftarrow}{\Delta}} X^{\cdot}$, and similarly that $Tot^n(X^{\mbox {\large $\cdot$}}) \cong\mbox{\hspace{.01cm}}\column{holim}{\column{\leftarrow}{\Delta^{(n)}}} X^{\cdot}$, where for any non-negative integer, $\Delta ^{(n)}$ denotes the full subcategory on the subsets of cardinality $\leq n$.  Also, let ${\cal D}_n$ denote the partially ordered set of non-empty subsets of the set $\underline{n}$, regarded as a category with a unique morphism from $S$ to $T$ whenever $ S \subseteq T$, and so that $Hom_{{\cal D}_n}(S,T) = \emptyset$ whenever $S \nsubseteq T$.  For any subset $S \subseteq \underline{n}$ of cardinality $s +1 $, let $\pi _n (S) = \underline{s}$.  Since $S$ inherits a total ordering from that of $\underline{n}$, there is a unique order-preserving bijective map $\xi _S : S \rightarrow \pi _n (S)$.   For any inclusion $S \subseteq T$ in ${\cal D}_n$, we  let $\pi _n  (S \subseteq T)$ denote the unique morphism in $\Delta ^{(n)}$ which makes the diagram

$$
\begin{diagram}\node{S}\arrow{e,t}{\xi _S} \arrow{s,t} {S \subseteq T}   \node{\pi _n  (S)} \arrow{s,b} {\pi _n  ( S \subseteq T)}\\
\node{T} \arrow{e,t}{\xi  _T} \node{\pi _n (T)}
\end{diagram}
$$

\noindent commute.  The two definitions  make $\pi _n $ into a functor from ${\cal D}_n$ to $\Delta ^{(n)}$.   Our next goal is now to prove that the natural pullback map 

$$ \pi _n ^* \mbox{\hspace{-.15cm}}: \mbox{\hspace{-.15cm}} \column{holim}{\column{\leftarrow}{\Delta ^{(n)}} } X^{\cdot} \rightarrow \column{holim}{\column{\leftarrow}{{\cal D}_n}}X^{\cdot} \raisebox{.05cm}{\tiny  {$\circ$ }}\pi _n 
$$

\noindent is an equivalence.  In order to do this,  we  recall that  by \cite{Bousfield}, Theorem 9.2, it will suffice show that for any object $\underline{k} \in \Delta ^{(n)}$, the category $\pi _n \downarrow \underline{k}$ has contractible nerve.  Recall that $ \pi _n \downarrow \underline{k}$ denotes the category whose objects are pairs $(S, \theta)$, where $S \in {\cal D}_n$ and $\theta : \pi _n (S) \rightarrow \underline{k}$ is a morphism in $\Delta ^{(n)}$, and where a morphism from $(S, \theta )$ to $(S^{\prime}, \theta ^{\prime})$  in $\pi _n \downarrow \underline{k}$ is a morphism $\varphi$ from $S$ to $S^{\prime}$ in ${\cal D}_n$ making the diagram 

$$
\begin{diagram} \node{\pi _n  (S) } \arrow{e,t}{\pi _n ( \varphi  )} \arrow{se,b}{\theta}  \node{\pi  _n (S^{\prime})}  \arrow{s,b}{\theta ^{\prime}}  \\
\node{}  \node{\underline{k}}
\end{diagram}
$$

\noindent commute.  We will construct a category equivalent to $\pi _n  \downarrow \underline{k}$ which is readily analyzed.  For any finite totally ordered set $S$, and any  positive integer $j$, a {\em $j$-fold interval decomposition} of  $S$ is a partition 

$$  S = S_0  \mbox{\hspace{-.1cm}\vspace{.1cm } \small $\coprod$}  S_1  \mbox{\hspace{-.1cm}\vspace{.1cm } \small $\coprod$}  \cdots  \mbox{\hspace{-.1cm}\vspace{.1cm } \small $\coprod$} S_ j
$$

\noindent of $S$,  so that whenever $0 \leq i < i^{\prime} \leq j$, then for any $x \in S_i$ and $y \in S_{i^{\prime}}$, $x < y$.   Note that some of the sets $S_i$ may be empty. We let $P_n^k$ be the partially ordered set whose objects are pairs $(S, \{S_j \}_{0 \leq j \leq k })$, where $S$ is a non-empty subset of $\underline{n}$, $\{ S_j \}_{0 \leq j \leq k}$ is a $k$-fold interval decomposition of $S$, and where 
$(S, \{S_j \}_{0 \leq j \leq k }) \leq (T, \{T_j \}_{0 \leq j \leq k })$ if and only if $S_j \subseteq T_j$ for all $0 \leq j \leq k$. 

\begin{Proposition} There is a natural equivalence of categories from the category attached to the poset $P_n^k$ to the category $\pi _n  \downarrow  \underline{k}$.  Consequently, in order to verify that $\pi _n  \downarrow  \underline{k}$ has contractible nerve, it suffices to show that $P_n^k$ does. 
\end{Proposition}
\begin{Proof}  Define $\alpha : P_n^k \rightarrow \pi _n  \downarrow \underline{k}$ by setting $\alpha ( (S, \{S_j \}_{0 \leq j \leq k })$ equal to the pair 
$(S, \theta )$, where $\theta :S \rightarrow \underline{k}$ is the unique order preserving map described by $\theta (x) = j$ if and only if $x \in S_j$.   It is easy to check that $\pi _n \downarrow \underline{k}$ is a partially ordered set in the sense that for any pair of objects $x$ and $y$, $Hom_{\pi _n \downarrow \underline{k}}(x,y)$ is either empty or consists of a single element.  One easily sees that the construction $\alpha$ respects the partial orderings,  and  hence creates a functor.  This process is clearly completely reversible (in fact, it is an {\em isomorphism} of categories), and the result follows.
\end{Proof} 

\begin{Proposition} $N_{\cdot}P_n^k$ is contractible whenever $k \leq n$.  
 \end{Proposition}
 
 \begin{Proof} We proceed by induction on $k$.  Of course, if $k=0$, we have that $P^0_n \cong {\cal D}_n$, and the latter poset has a maximal element, whence its nerve is contractible.  Now, suppose we know the result to be true for all $k \leq K$, and we wish to prove it for $K$.   Let $n \geq K$.  The poset $P^K_n$ contains $K+1$ minimal elements $\{ x_0, x_1, \ldots , x_K \}$, where $x_j$ denotes the element whose underlying set is $\{ n \}$,  and whose partition places $n$ in the $j$-th interval.  Let $\tilde{P}^K_n \subseteq P^K_n$ denote the partially ordered subset 
 $P^K_n - \{ x_0, x_1, \ldots , x_K \}$.  We also have the obvious embedding $i : P_{n-1}^K \hookrightarrow \tilde{P}_n^K$. We claim that this inclusion induces a homotopy equivalence on nerves.  To see this, it suffices to construct an order preserving map $f : \tilde{P}^K_n \rightarrow P_{n-1}^K$ so that $f \compcirc i = id$, and so that $i \compcirc f (x) \leq x$ for all $x \in \tilde{P}_n^k$.  But now we note that the map given by $(S, \{ S_j \}) \rightarrow (S- \{n \}, \{ S_j - \{n \} \})$ provides the required map of partially ordered sets.  For each $0 \leq j \leq K$, we let $Q_j \subseteq \tilde{P}_n^K$ denote the subset $\{ y \in \tilde{P}_n^K \mid y \geq x_j \}$. As in \cite{Quillen4}, one now verifies that we have a decomposition 
 
 $$ N_{\cdot} P_n^K \cong N_{\cdot} \tilde{P}_n^K  \cup C N_{\cdot} Q_0  \cup CN_{\cdot} Q_1 \cup \cdots  \cup C N_{\cdot}Q_K
 $$
 
 \noindent We next observe that the inclusion $Q_0 \hookrightarrow \tilde{P}_n^K$ induces an equivalence on nerves.   To see this, we observe that the restriction of the poset map $f$ constructed above to $Q_0$ is an isomorphism of partially ordered sets. It follows that the simplicial set $N_{\cdot} \tilde{P}_n^K  \cup C N_{\cdot} Q_0$ is contractible.  Now, $N_{\cdot} P_n^K$ is obtained from  $N_{\cdot} \tilde{P}_n^K  \cup C N_{\cdot} Q_0$ by adjoining cones to the subspaces $N_{\cdot}Q_1, N_{\cdot} Q_2 , \ldots ,  N_{\cdot} Q_n$.   This means that in order to prove the contractibility of $N_{\cdot} P_n^K$, it will suffice to prove the contractiblity of the sets $N_{\cdot} Q_j$, for $j = 1, 2, \ldots , n$.  We claim that the partially ordered set $Q_j$ is isomorphic to $P_{n-1} ^{K-j}$.  Since we have $K-j \leq n-1$ whenever $j \geq 1$,  we have that $N_{\cdot} P_{n-1}^{K-j}$ is contractible by the inductive hypothesis, and the result would follow.   To establish the isomorphism 
 $Q_j \cong P_{n-1}^K$, we first note that $(S, \{ S_i \}_{0 \leq i \leq K } ) \in Q_j$ if and only if $n \in S_{K-j}$.  Therefore, for $(S, \{ S_i \}_{0 \leq i \leq K } ) \in Q_j$ we have that $S_l = \emptyset $ for $l > K-j$.  The isomorphism of partially ordered sets is now given by $(S, \{ S_i \}_{0 \leq i \leq K } )  \rightarrow (S -\{ n \}, \{ S_{i} - \{ n \} \} _{0 \leq i \leq K-j})$
  \end{Proof}
  
  \noindent Consider any  spectrum  valued functor $F$ on  ${\cal D}_n$, and  for any spectrum $X$,  we let $F_X$ denote the constant functor on ${\cal D}_n$ with value $X$.  Suppose that we are given a natural transformation $F_X \rightarrow F$.  We are interested in constructing a useful model for the homotopy fiber of the natural map 
 
  $$ \holim{{\cal D}_n}{F_X} \rightarrow \holim{{\cal  D}_n}{F}
  $$
  
  \noindent We let $I$ denote the partially ordered set of non-empty subsets of $\underline{1}$, i.e. the category pictured by the diagram 
  
  $$
  \begin{diagram} \node{ \{ 0 \} } \arrow{e} \node{\{0,1 \}}\node {\{ 1 \} }\arrow{w}
  \end{diagram}
  $$
  \noindent and let $J$ denote the full subcategory on the two objects $\{ 0 \} $ and  $\{ 0,1 \}$. 
Consider the category $I^{n+1}$,  and its full subcategory $J^{n+1}$.   For  any object $\sigma = ( S_0 , S_1 ,\cdots , S_{n})$ of $J^{n+1}$,  we let $\psi (\sigma) \subseteq \underline{n} $ denote the subset $\{ j \in \underline{n} \mid S_j = \{ 0,1 \} \}$.  We now define a spectrum valued functor $\overline{F}$ on $I^{n+1}$ by the formulas
\renewcommand{\arraystretch}{1.5}

$$
\left\{ \begin{array}{l}
\overline{F}(0, 0, \ldots, 0) = X \\
\overline{F}( \sigma ) = F(\psi (\sigma))\mbox{ for any } \sigma \in J^{n+1} - (0, 0, \ldots, 0 ) \\
\overline{F} (\xi ) = * \mbox{ for any } \xi \in I^{n+1} - J^{n+1} 
\end{array}
\right.
$$

\noindent The behavior on morphisms is evident, when we recall that we are given a natural transformation from the constant functor $F_X$  to $F$.  

\begin{Proposition} \label{fiber}There is a natural equivalence (in the homotopy category)  $\varepsilon _n$ from the homotopy fiber $\Phi$ of  the natural map 
$\holim{{\cal D}_n}{F_X} \rightarrow \holim{{\cal  D}_n}{F}$ to $\holim{{I^{n+1}}}{\overline{F}}$.  %Moreover, there is a homotopy commutative diagram 

%$$
%\begin{diagram} \node{Fiber(\holim{{\cal D}_{n+1}}{F_X} \rightarrow \holim{{\cal  D}_{n+1}}{F})}
%\arrow{e,t}{\varepsilon _{n+1}}\arrow{s}  \node{\holim{{I^{n+2}}}{\overline{F}}} \arrow{s}  \\
%\node{Fiber(\holim{{\cal D}_{n}}{F_X} \rightarrow \holim{{\cal  D}_{n}}{F})} \arrow{e,t} {\varepsilon _{n}} 
%\node{\holim{{I^{n+1}}}{\overline{F}}}
%\end{diagram}
%$$

\end{Proposition}
\begin{Proof}  $\Phi$ can clearly be interpreted as the homotopy inverse  limit over the category ${\cal D}_n \times I$ of the functor $G$ defined by 
$$
\left\{ \begin{array}{l}
G(S,\{ 0 \}) = X  \\
G(S, \{ 0,1 \}) = F(S) \\
G(S, \{ 1 \} ) = *
\end{array}
\right.
$$

Define ${\cal E}_n$ to be the partially ordered set obtained from the product poset ${\cal D}_n \times  I$ by identifying the subset $\{ 0 \} \times {\cal D}_n$ to a single point $\epsilon$.  This means that the object set of ${\cal E}_n$ is $\{ \epsilon \} \cup {\cal D}_n \times \{ 0,1 \} \cup {\cal D}_n \times \{ 1 \}$, that $\epsilon \leq (S, \{ 0,1 \} )$ for all $S \in {\cal D}_n$, that $\epsilon$ is incomparable with any element of the form $( S, \{ 1 \})$, and that the ordering relationships between any elements of ${\cal D}_n \times \{ \{0,1\}, \{ 1 \} \}$ are identical to those in the original set ${\cal D}_n \times I$.  There is a natural projection $\pi : {\cal D}_n \times I \rightarrow {\cal E}_n$.  It is clear from the definition that the functor $G$ factors over $\pi$, i.e. that there is a spectrum valued functor $\overline{G}$ on ${\cal  E}_n$ so that $G = \overline{G} \raisebox{.03 cm}{\tiny $\circ$} \pi$.   On the other hand, we may also define a functor $\rho : I^{n+1} \rightarrow {\cal E}_n$ as follows.  For any element $v = \{ S_0, S_1, \ldots , S_n \}$ of $I^{n+1} - J^{n+1}$, we define $\theta (v) \in {\cal D}_n$ to be 
the subset $\{ j \in \underline{n} \mid S_j = \{ 1 \} \mbox{ or } \{ 0,1 \} \}$.  The functor $\rho$ is now defined by 
$$
\left\{ \begin{array}{l}
\rho ( \{ 0 \}, \{ 0 \}, \ldots , \{ 0 \} ) = \epsilon \\
\rho ( S_0, S_1, \ldots S_n ) = \psi (S_0, S_1, \ldots , S_n) \times \{ 0, 1 \} \mbox{ for any } (S_0, S_1, \
\ldots, S_n ) \in J^{n+1}- (\emptyset, \emptyset, \ldots , \emptyset ) \\
\rho ( S_0, S_1, \ldots , S_n) = \theta (S_0, S_1, \ldots , S_n ) \times \{ 1 \} \mbox{ for any } 
(S_0, S_1, \ldots , S_n ) \in I^{n+1} - J^{n+1}
\end{array}
\right.
$$ 
Behavior on morphisms is determined since all categories in question are partially ordered sets, and one readily checks that the ordering is respected. One now checks directly that $G \raisebox{.03cm}{\tiny $\circ$} \rho$ is identical to the functor $\overline{F}$ defined above. It is also direct to check that for any object $ x \in {\cal E}_n$, the categories $\pi \downarrow x$ and $\rho \downarrow x$ have contractible nerves.  The result now follows from \cite{Bousfield}, Theorem 9.2.\end{Proof}

\noindent Now let $A$ be a commutative  $S$-algebra, $B$ a commutative $A$-algebra spectrum,  and let $\Psi$ denote the $I$-diagram

$$  A \longrightarrow B  \longleftarrow *
$$

\noindent  of $A$-module spectra, i.e $\{ 0 \} \rightarrow A$, $\{ 0,1 \} \rightarrow B$, and $\{ 1 \} \rightarrow *$.  Then we can define an $I^{n+1}$-diagram $\Psi^n$ by 

$$\Psi^n(S_0,S_1, \ldots , S_n) =  \Psi (S_0) \column{\wedge}{A} \Psi (S_1) \smashover{A} \ldots \smashover{A} \Psi (S_n )
$$
This is a diagram of $A$-bimodule spectra.  We record the following lemma concerning homotopy inverse limits of these $I$-diagrams.  

\begin{Lemma} \label{smash1} Let $M$ denote any  $A$-module.  Then 
$$ \holim{I}{(\Psi(-) \smashover{A} M)} \cong (\holim{I}{\Psi}) \smashover{A} M
$$
\end{Lemma}
\begin{Proof}  We note that by the definitions of homotopy inverse limits, we find that 
$\holim{I}{\Psi}$ is the homotopy fiber of the spectrum map $A \rightarrow B$.  Consequently, the result reduces to  the fact that given any map $f : P \rightarrow Q$ of  $A$-module spectra, and any $A$-module spectrum $M$, we have an equivalence $ \mbox{Fib}(f) \smashover{A} M \cong \mbox{Fib}(f \smashover{A} id_M)$.  But this is clear from the results of  \cite{May et all} and \cite{Smith et all}. 
\end{Proof}

\noindent The following result will be key in proving our main theorem. 

\begin{Proposition}\label{smash2} Let $A$ and $\Psi$ be as above, and let $M$ be a $A$-module spectrum.  .  Then we have a natural equivalence
$$ \holim{{I^{n+1}}}{M \smashover{A} \Psi^n} \cong M \smashover{A} \underbrace{\holim{I}{\Psi} \smashover{A} \cdots \smashover{A} \holim{I}{\Psi}}_{n+1 \mbox{ factors }}
$$

\end{Proposition}
\begin{Proof} By induction on $n$.  The result is trivially true for $n = 0$.  Suppose that it holds for $n = N$, and we wish to prove it for $n = N+1$.  We have the string of equivalences

$$ \begin{array}{l} \holim{\mbox{\tiny $(S_0, \ldots S_{N+1}) \in I^{N+2}$}}{(M \smashover{A} \Psi^{N+2}(S_0, \ldots , S_{N+1}))}  \\
\cong \mbox{ } \holim{S_0 \in I}{\left[ \holim{(S_1, \ldots , S_{N+1}) \in I^{N+1}}{M \smashover{A} \Psi (S_0 ) \smashover{A} \Psi ^{N+1}(S_1, \ldots , S_{N+2})}\right]}  \\
\cong \mbox{ } \holim{S_0 \in I}{M \smashover{A}\Psi (S_0 ) \smashover{A} \left[  \underbrace{\holim{I}{\Psi} \smashover{A} \cdots \smashover{A} \holim{I}{\Psi}}_{N+1 \mbox{ factors }} \right] 
} \\
\cong \mbox{ }M \smashover{A}\left[ \holim{I} {\Psi}\right] \smashover{A}\left[   \underbrace{\holim{I}{\Psi} \smashover{A} \cdots \smashover{A} \holim{I}{\Psi}}_{N+1 \mbox{ factors }} \right] 
\\
\cong \mbox{ }M \smashover{A} \left[   \underbrace{\holim{I}{\Psi} \smashover{A} \cdots \smashover{A} \holim{I}{\Psi}}_{N+2 \mbox{ factors }} \right]

\end{array}
$$
\noindent The first equivalence is a ``Fubini'' type theorem for homotopy inverse limits, the second is an application of the inductive hypothesis, and the third is an application of \ref{smash1}. 
\end{Proof}

\begin{Corollary} \label{interpretation} Let $f : A \rightarrow B$ be a homomorphism of commutative  $S$-algebras,  and let $M$ be a  $A$-module spectrum.  Let ${\cal I}$ denote the homotopy fiber of the homomorphism $f$.  ${\cal I}$ is an $A$-module spectrum.  Then there is a homotopy fibration sequence of $A$-module spectra 
$$ {\cal I}^{\raisebox{.05cm}{ \hspace{-.1cm}\tiny $\smashover{A}$} (k+1)} \smashover{A}M \rightarrow M  \rightarrow Tot^k T^{\cdot}_A(M;B)$$
Moreover, these sequences fit together in  homotopy commutative diagrams 

\begin{equation} \label{consistent}
\begin{diagram}
\node{ {\cal I}^{\raisebox{.05cm}{ \tiny $\smashover{A}$} (k+2)} \smashover{A} M} \arrow{e} \arrow{s} \node{M} \arrow{e} \arrow{s,t} {id_M}  \node{Tot^{k+1}(T^{\cdot}_A(M;B)} \arrow{s} \\
\node{ {\cal I}^{\raisebox{.05cm}{ \hspace{-.1cm}\tiny $\smashover{A}$} (k+1)} \smashover{A} M} \arrow{e} \node{M} \arrow{e} \node{Tot^k(T^{\cdot}_A(M;B)}
\end{diagram}
\end{equation}
where the right hand vertical map is the usual projection, and the left hand vertical map is the composite

\begin{equation} \label{explain}
  {\cal I}^{\raisebox{.05cm}{ \hspace{-.1cm}\tiny $\smashover{A}$} (k+2)}  \cong {\cal I} \smashover{A}  {\cal I}^{\raisebox{.05cm}{\hspace{-.1cm} \tiny $\smashover{A}$} (k+1)} \stackrel{j \smashover{A} id}{\longrightarrow}   A \smashover{A} {\cal I}^{\raisebox{.05cm}{ \hspace{-.1cm}\tiny $\smashover{A}$} (k+1)} \cong  {\cal I}^{\raisebox{.05cm}{ \hspace{-.1cm}\tiny $\smashover{A}$} (k+1)} 
\end{equation}
$j:{\cal I} \rightarrow A$ being the evident inclusion.

\end{Corollary}

This result implies the following technical corollary.

\begin{Corollary} \label{connectivity}  Suppose that we have a homomorphism $f:A \rightarrow B$ of $(-1)$-connected commutative  $S$-algebras, and that $\pi _0 (f)$ is surjective.  Suppose that $M$ is a $k$-connected $A$-module spectrum.  Then each of the spectra $Tot^n({\cal T}^{\cdot}_A(M;B))$ is $k$-connected, and  $M_B = Tot({\cal T} ^{\cdot}_A(M;B))$ is at least $k-1$-connected. 
\end{Corollary}
\begin{Proof}  We prove the first part by induction on $n$.  The result holds for $n=0$ since $Tot^0({\cal T}^{\cdot}_A(M;B)) \cong B \smashover{A}M$.  The K\"{u}nneth spectral sequences for smash products shows that the connectivity $B \smashover{A}M$ is at least $k$. For higher values of $n$, we suppose the result know for the value $n-1$.  By the commutative diagram \ref{consistent} in 
Corollary \ref{interpretation} above, we have an identification of the homotopy fiber of the map 
$Tot ^n ({\cal T}^{\cdot}_A(M;B)) \rightarrow Tot ^{n-1} ({\cal T}^{\cdot}_A(M;B))$ with the homotopy cofiber of the   map $I^{\raisebox{.05cm}{ \hspace{-.1cm}\tiny $\smashover{A}$} (n+1)} \rightarrow I^{\raisebox{.05cm}{ \hspace{-.1cm}\tiny $\smashover{A}$} n}$  in diagram \ref{explain} above.  Since both $I^{\raisebox{.05cm}{ \hspace{-.1cm}\tiny $\smashover{A}$} (n+1)} $ and $I^{\raisebox{.05cm}{ \hspace{-.1cm}\tiny $\smashover{A}$} n}$ are (-1)-connected, so is the homotopy fiber of the map $Tot ^n ({\cal T}^{\cdot}_A(M;B)) \rightarrow Tot ^{n-1} ({\cal T}^{\cdot}_A(M;B))$.  It now readily follows that that $Tot^n({\cal T}^{\cdot}_A(M;B))$ is (-1)-connected.  The result for $Tot({\cal T}^{\cdot}_A(M;B))$ now follows from the $lim^1$-exact sequence for homotopy fibrations (see \cite{Bousfield}, Ch X, Section 3).
\end{Proof}

\noindent We will also need the following technical lemma.

\begin{Lemma} \label{split} Let $A \rightarrow B \stackrel{f}{\rightarrow} C$ be a diagram of commutative  $S$-algebras, and let $M$ denote an  $A$-module spectrum.  Suppose further that there is a homomorphism of commutative $A$-algebra spectra $s: C \rightarrow B$ such that $s \compcirc f = id_B$.  Then the natural map $M_B^{\wedge} \rightarrow M_C^{\wedge}$ is an equivalence of spectra. 
\end{Lemma}
\begin{Proof} Follows directly from Proposition \ref{BarrBeckcorollary}, applied to the triples $S = B \smashover{A}-$ and $C \smashover{A}-$.  
\end{Proof}

{\bf Proof (of Theorem \ref{main}):} We consider first the case when the diagram of $S$-algebra  homomorphisms is $A \rightarrow A \rightarrow C$, i.e. $A = B$, and where the homotopy fiber ${\cal I}$ of the map $f: A= B \rightarrow C$ is 0-connected, i.e. when $\pi _0 (f) $ is an isomorphism and $\pi _1 (f)$ is surjective.  It follows from Propositions \ref{fiber} and \ref{smash2} that the homotopy fiber ${\cal F}$ of the map $ M \rightarrow \holim{{\Delta^{n}}} {{\cal T}^{\cdot}_A(M;B) }\cong Tot^n {\cal T}^{\cdot}_A(M;B)$ is equivalent to $\underbrace{ {\cal I} \smashover{A} \cdots \smashover{A}  {\cal I} }_{n+1 \mbox{ factors }}$. Since ${\cal I}$ is 0-connected, it follows by an easy applicaton of Proposition \ref{oneconnectivity} that 
${\cal F}$ is $n$-connected. This gives the result in this case, since when one passes to the homotopy limit over $n$ the connectivity of the map goes to infinity, and since $M^{\wedge}_A \cong M$. 

\noindent We next consider the case where $A$ is not necessarily equal to $B$, but where we have that the homotopy fiber ${\cal I}$ of the map $f:B \rightarrow C$ is 0-connected.   We construct the bicosimplicial  spectrum ${\cal C} ^{\cdot \cdot}$ given by ${\cal C}^{pq} = {\cal T}^{p}_A( {\cal T}^{q}_A(M;B);C)$, with the obvious coface and codegeneracy maps.  Lemma \ref{insertion} asserts that   the natural map 
$${\cal  T}^{\cdot}_A( M;C) \rightarrow {\cal C}^{\cdot \cdot}
$$
(induced by $\eta: M \rightarrow {\cal T}_A^{\cdot}(M;B)$)
induces an equivalence on total spectra, where ${\cal T}^{\cdot}_A( M;C)$ is regarded as a bicosimplicial spectrum constant in the $q$-direction. We also have the natural map $\eta: {\cal T}^{\cdot}_A(M;B) \rightarrow  {\cal C}^{\cdot \cdot}$ obtained by regarding ${\cal T}^{\cdot}_A(M;B)$ as a bicosimplicial spectrum constant in the $p$-direction.  Since $Tot ({\cal C}^{\cdot \cdot}) \cong M^{\wedge}_C$, and this equivalence is compatible with the standard map $M_B^{\wedge} \rightarrow M_C^{\wedge} $, it suffices to show that $\eta $ induces a weak equivalence on total spectra.  For this it suffices to show that each of the maps 
$$ {\cal T}^{q}_A(M;B) \rightarrow Tot({\cal T} ^{\cdot}_A({\cal T} ^{q}_A(M;B));C)
$$
is a weak equivalence of spectra.  Since each of the $A$-module spectra 
$${\cal T} ^{q}_A(M;B) =  \underbrace{B \smashover{A} \cdots \smashover{A} B}_{q+1 \mbox{ factors }} \smashover{A} M
$$
is equipped with a $B$-module structure restricting to its given $A$-module structure, it will now suffice to prove that the natural map $M_B^{\wedge} \rightarrow M_C^{\wedge} $ is an equivalence for any $B$-module spectrum $M$.  For a $B$-module spectrum $M$, we may construct the cosimplicial spectrum ${\cal T}^{\cdot}_B (M ; C)$.  From the hypothesis on $f$, and by the previous case, we find that the connectivity of the natural maps $M \rightarrow Tot^n({\cal T}^{\cdot}_B(M;C))$ goes to infinity with $n$. We will now apply the functors ${\cal T} ^{\cdot}_A(-;B)$ and ${\cal T} ^{\cdot}_A (-;C)$ levelwise to the cosimplicial spectrum ${\cal T}^{\cdot}_B (M;C)$ to obtain bicosimplicial spectra, with a natural bicosimplicial map  from the first to the second induced by the $A$-algebra  spectrum map $B \rightarrow C$.   We obtain a commutative diagram 

\begin{equation}\label{compcompare}
\begin{diagram} \node{M_B^{\wedge}  = \mbox{\hspace{.01cm}} \column{Tot}{p \in \Delta} {\cal T}^p_A(M;B)} \arrow{s} \arrow{e} \node{M_C^{\wedge} =\mbox{\hspace{.01cm}} \column{Tot}{p \in \Delta} {\cal T}^p_A(M;C)} \arrow{s}\\
\node{\column{Tot}{p \in \Delta} \mbox{\hspace{.01cm}}\column{Tot}{q \in \Delta} {\cal T} ^p_A({\cal T}^q_B(M;C);B)} \arrow{e}
\node{\column{Tot}{p \in \Delta}  \mbox{\hspace{.01cm}}\column{Tot}{q \in \Delta} {\cal T}^p_A({\cal T} ^q_B(M;C);C)}
\end{diagram}
\end{equation}
We wish to prove that the upper horizontal arrow is an equivalence of spectra.  We will do this by showing that the two vertical arrows and the lower horizontal arrow are equivalences of spectra. We first note that by the increasing connectivity (with $n$)  of the maps $M \rightarrow \column{Tot}{q \in \Delta}^n {\cal T}^q_B(M;C)$ and Corollary \ref{connectivity}, we find that the natural maps 
$$ \column{Tot}{p \in \Delta} {\cal T}^p_A(M;B) \rightarrow \holim{n}{\column{Tot}{p \in \Delta} {\cal T} ^p_A( \column{Tot}{q \in \Delta} ^n {\cal T} ^q _B (M;C);B)}
$$ 
and
$$ \column{Tot}{p \in \Delta} {\cal T}^p_A(M;C) \rightarrow 
\holim{n}{\column{Tot}{p \in \Delta} {\cal T} ^p_A( \column{Tot}{q \in \Delta^n }{\cal T} ^q _B (M;C);C)}
$$  
are equivalences of spectra. 
By Proposition \ref{smashproperty}, we have that $${\cal T} ^p_A( \column{Tot}{q \in \Delta} ^n {\cal T}^q _B (M;C);B) \cong  \column{Tot}{q \in \Delta} ^n {\cal T} ^p_A({\cal T} ^q _B (M;C);B)$$ and $${\cal T}^p_A( \column{Tot}{q \in \Delta} ^n {\cal T} ^q _B (M;C);C) \cong  \column{Tot}{q \in \Delta} ^n {\cal T} ^p_A({\cal T}^q _B (M;C);C)$$
Therefore, the maps 
$$ \column{Tot}{p \in \Delta} {\cal T}^p_A(M;B) \rightarrow \holim{n}{\column{Tot}{p \in \Delta}  \column{Tot}{q \in \Delta} ^n {\cal T}^p_A({\cal T} ^q _B (M;C);B)} \cong \column{Tot}{p \in \Delta} \mbox{\hspace{.01cm}}\column{Tot}{q \in \Delta} {\cal T} ^p_A({\cal T} ^q_B(M;C);B)
$$ 
and
$$ \column{Tot}{p \in \Delta} {\cal T} ^p_A(M;C) \rightarrow \holim{n}{\column{Tot}{p \in \Delta}  \column{Tot}{q \in \Delta} ^n {\cal T} ^p_A({\cal T}^q _B (M;C);C)} \cong \column{Tot}{p \in \Delta} \mbox{\hspace{.01cm}}\column{Tot}{q \in \Delta} {\cal T} ^p_A({\cal T} ^q_B(M;C);C)
$$ 
are also weak equivalences of spectra.  These are the vertical arrows in Diagram \ref{compcompare}.  We must also check that the lower horizontal arrow in Diagram \ref{compcompare} is an equivalence of spectra.  This arrow is of course induced by a map of bicosimplicial spaces, so we may verify it by proving that it is a levelwise equivalence.  Specifically, if we can show that for each $q$, the  map 
$$ \column{Tot}{p \in \Delta}{\cal T} ^p_A({\cal T} ^q_B(M;C);B) \rightarrow
\column{Tot}{p \in \Delta}{\cal T} ^p_A({\cal T} ^q_B(M;C);C) 
$$
is an equivalence of spectra, then the lower horizontal map in Diagram \ref{compcompare} will be an equivalence.  But, ${\cal T} ^q_B(M;C)$ admits a $C$-module structure (and therefore also a $B$-module structure), so by Proposition \ref{properties},(6), we have
$$  \column{Tot}{p \in \Delta}{\cal T} ^p_A({\cal T} ^q_B(M;C);B) \cong \column{Tot}{p \in \Delta}{\cal T} ^p_A({\cal T}^q_B(M;C);C) 
$$
which gives the result, i.e. that the statement of the theorem holds in the case when ${\cal I}$ is 0-connected.  

\noindent The general case, i.e. without the assumption that $\pi _1 B \rightarrow \pi _1 C$ is surjective, can now be handled as follows.  Let $\{ f_{\alpha} \}_{\alpha \in a}$ be a generating set for $\pi _1 C$.  Let $Z = \bigvee _{\alpha \in A} S^1_{\alpha}$denote the suspension spectrum of a bouquet of circles parametrized by $A$, and let $\theta : Z \rightarrow C$ denote the map whose restriction to  $S^1_{\alpha}$ is the map $f_{\alpha}$. Let $B(Z)$ denote the free $B$-module on the spectrum $Z$.  Then we obtain a map of commutative $A$-algebra spectra $g: Sym_B ( B(Z)) \rightarrow C$. One readily checks that $\pi_0(g)$ is an isomorphism and that $\pi _1 (g)$ is surjective, so we may conclude that $M_{Sym_B(B(Z))} \rightarrow M_C$ is an equivalence of spectra.  On the other hand, the inclusion $B \rightarrow Sym_B (B(Z))$ admits a section, and so the natural map $M_B \rightarrow M_{Sym_B(B(Z))}$ is also an equivalence by \ref{split}.  This concludes the proof of Theorem \ref{main}.\qed

\noindent We have been discussing the invariance of completion along a $S$-algebra homomorphism $f : A \rightarrow B$ under changes in the target $S$-algebra $B$.  We will also need a similar result which will demonstrate invariance under changes in $A$ under suitable circumstances.  For a homomorphism of commutative  $S$-algebras $f : A \rightarrow B$, and a $B$-module spectrum $M$, we will denote by $\rho ^AM$ the result of regarding $M$ as an $A$-module spectrum via the $S$-algebra  homomorphism $f$.  

\begin{Theorem} \label{ringinvar} Let $A \stackrel{f}{\rightarrow} B \rightarrow C$ be a diagram of $(-1)$-connected commutative  $S$-algebras, and let $M$ be a $B$-module spectrum.  Suppose that the homomorphism $\pi _0 f : \pi _0 A \rightarrow \pi _0 B$ is an isomorphism.  Then the natural homomorphism $(\rho ^A M)_C^{\wedge} \rightarrow M_C^{\wedge}$ is a weak equivalence of spectra.   
\end{Theorem}
\begin{Proof}  Let $T_A$ and $T_B$ denote the triples $M \mapsto C \smashover{A} M$ and $M \mapsto C \smashover{B} M$, respectively.  We consider the bicosimplicial $A$-module spectrum 
$$ (p,q) \mapsto T_B^p T_A^q (M)
$$
as usual, and we have by Lemma \ref{insertion} that its total space is equivalent to $ M^{\wedge}_C$.  In order to check the equivalence we are interested in, it will be sufficient to prove that for each $q$, the total spectrum of the cosimplicial spectrum
$$p \mapsto T_B^p T_A^q (M)
$$
is weakly equivalent to $T_A^q (M)$ under the natural map $\eta _B (T_A^q (M))$.   As in the proof of Theorem \ref{main}, we let 
${\cal I} $ denote the  homotopy fiber of the homomorphism $B \rightarrow C$.  ${\cal I}$ is a $B$-module.  Then Corollary \ref{interpretation} shows that it will suffice to prove that the inverse system of spectra 
$$  \cdots \rightarrow {\cal I}^{\smashover{B}^{{k+1}}} \smashover{B} T^p _A(M)\rightarrow 
{\cal I}^{\smashover{B}^{{k}}}  \smashover{B} T^p _A(M)\rightarrow {\cal I}^{\smashover{B}^{{k-1}}}  \smashover{B} T^p _A(M) \rightarrow \cdots 
$$
has contractible homotopy inverse limit.  It is clear that to verify this, it will suffice to show that for any $B$-module spectrum $N$, the inverse system of spectra 
$$  \cdots \rightarrow{\cal I}^{\smashover{B}^{{k+1}}} \smashover{B} C \smashover{A}N\rightarrow 
{\cal I}^{\smashover{B}^{{k}}}  \smashover{B} C \smashover{A}N\rightarrow {\cal I}^{\smashover{B}^{{k-1}}}  \smashover{B} C \smashover{A}N \rightarrow \cdots 
$$
has contractible homotopy inverse limit.  In order to prove this, we observe that there is a standard equivalence
$$ C \smashover{A} N \cong C \smashover{B} B \smashover{A} B \smashover{B} M
$$
of $B$-module spectra.  The idea will be to replace $B \smashover{A} B$ by its cosimplicial approximation ${\cal T}^{\cdot}_{B \smashover{A}B}(B \smashover{A}B; B)$, where $B$ is a commutative $B \smashover{A}B$-algebra via the multiplication map $B \smashover{A}B \rightarrow B$.  Since we are assuming that $\pi _0f$ is an isomorphism, it follows that the map $\pi _0 B \smashover{A} B \rightarrow \pi _0 B$ is an isomorphism, and therefore by Theorem \ref{main} that the natural map $\eta$ gives an equivalence $B \smashover{A} B \rightarrow (B \smashover{A} B )^{\wedge}_{B}$.  Also, $\pi _1 B \smashover{A}B \rightarrow \pi _1 B$ is surjective since the multiplication map has a section.  It follows that 
the natural map $ B \smashover{A}B \rightarrow Tot^i( {\cal T}^{\cdot}_{B \smashover{A} B}(B \smashover{A} B ; B))$ is at least $(i-1)$-connected.    Letting $T^i = Tot^i( {\cal T}^{\cdot}_{B \smashover{A} B}(B \smashover{A} B ; B))$, it now follows from the K\"{u}nneth spectral sequence that the natural map 
$C \smashover{A} N \rightarrow C \smashover{B} T^i \smashover{B} N$ is $(i-1)$-connected.  Consequently, we have an equivalence 
$$ C \smashover{A} N \cong Tot(k \mapsto C \smashover{B} {\cal T}^k _{B \smashover{A} B}(B \smashover{A} B; B) \smashover{B} N)
$$
Therefore, it will now suffice to prove that the homotopy inverse limit of the system 
$$  \cdots \rightarrow {\cal I}^{\smashover{B}^{{k+1}}} \smashover{B} C  \smashover{B} {\cal T}^k _{B \smashover{A} B}(B \smashover{A} B; B) \smashover{B}N\rightarrow 
{\cal I}^{\smashover{B}^{{k}}}  \smashover{B} C  \smashover{B} {\cal T}^k _{B \smashover{A} B}(B \smashover{A} B; B) \smashover{B}N\rightarrow $$
$$ \rightarrow {\cal I}^{\smashover{B}^{{k-1}}}  \smashover{B} C  \smashover{B} {\cal T}^k _{B \smashover{A} B}(B \smashover{A} B; B) \smashover{B}N \rightarrow \cdots 
$$
is contractible. 
Now, ${\cal T}^k _{B \smashover{A} B}(B \smashover{A} B; B)$ is a $(k+1)$-fold smash product over $ B \smashover{A}B$ of copies of $B$. By including $B$ into the first tensor factor, we find that it can be viewed as a $B$-algebra, for which we write $\Lambda$.  $\Lambda$ can be viewed as a $B-B$-bimodule, with both right and left actions coming from the (left) action of $B$ on $\Lambda$.  We can now resolve $\Lambda$ by free $B$ -modules (each of which is regarded as a bimodule via the left action, and which can be taken to be free on a wedge of spheres of the same dimension), and it is now easy to see that it will suffice to show that the inverse system 
$$  \cdots \rightarrow {\cal I}^{\smashover{B}^{{k+1}}} \smashover{B} C  \smashover{B} F \smashover{B}N\rightarrow 
{\cal I}^{\smashover{B}^{{k}}}  \smashover{B} C  \smashover{B} F \smashover{B}N 
 \rightarrow {\cal I}^{\smashover{B}^{{k-1}}}  \smashover{B} C  \smashover{B} F \smashover{B}N \rightarrow \cdots 
$$
has contractible homotopy inverse limit, where $F$ is any free $B$-module on a bouquet of spheres of a fixed dimension $k$, regarded as a $B - B$ bimodule, with both actions agreeing with the left action on $F$.  But to check this, it will suffice to show that the map 
$${\cal I}^{\smashover{B}^{{k}}}  \smashover{B} C  \smashover{B} F \smashover{B}N 
 \rightarrow {\cal I}^{\smashover{B}^{{k-1}}}  \smashover{B} C  \smashover{B} F \smashover{B}N  
$$
is null homotopic as a map of spectra. But this question clearly can be reduced to the case where $F = B$, so we will have to show that the map 
$${\cal I}^{\smashover{B}^{{k}}}  \smashover{B} C   \smashover{B}N 
 \rightarrow {\cal I}^{\smashover{B}^{{k-1}}}  \smashover{B} C   \smashover{B}N  
$$
is null homotopic.  This can clearly be reduced to the case where $k = 1$, so we need to  check that the map 
$${\cal I} \smashover{B} C \smashover{B}N \rightarrow 
B \smashover{B} C \smashover {B} N \cong C \smashover{B}N
$$
is null.  But we have the cofiber sequence 
$$ {\cal I} \smashover{B} C \smashover{B}N \rightarrow 
 C \smashover {B} N \rightarrow C \smashover{B} C \smashover{B} N
$$
and the second map is the  inclusion of a wedge summand, since we have the retraction 
$$ \mu \smashover{B} id_N: C \smashover{B} C \smashover{B} N \rightarrow  C \smashover {B} N
$$
where $\mu:  C \smashover{B} C \rightarrow C$ is the multiplication map.  This gives the result.  \end{Proof}

\section{Algebraic to geometric spectral sequence} \label{algeo}
In this section, we will provide a computational device which will permit the computation of homotopy groups of derived completions of module spectra over commutative  $S$-algebras in terms of derived completions over actual rings.  Here is the statement of the main theorem of this section. 

\begin{Theorem} \label{algebraictogeometric} Let $A \rightarrow B$ be a map of commutative (-1)-connected  $S$-algebras, with $\pi _0 A \rightarrow \pi _0 B$ surjective,  and let $M$ be an $A$-module spectrum.   $\pi _0A$ is a commutative ring, $\pi _0 B$ is a commutative $\pi _0A$-algebra, and for each $i$, $\pi _iM$ is a $\pi _0 A$-module.  We may therefore construct the derived completion spectrum $(\pi _iM)_{\pi _0 B}$ for each $i$.   There is a second quadrant spectral sequence with $E_1^{pq} = \pi _{q+2p} ((\pi _{-p} M)_{\pi _0B}^{\wedge})$, converging to $\pi _{p+q}(M_B)$.
\end{Theorem}

\noindent We will require some preliminary technical work on Postnikov decompositions of $S$-algebras and module spectra before we can present the proof of this result.   We recall from \cite{shipley} that the category  $\mbox{Alg}_S$ admits the structure of a Quillen model category, in which the weak equivalences are the homomorphisms of  $S$-algebras inducing isomorphisms on homotopy groups.  There is also a ``free commutative  $S$-algebra" on a spectrum $X$, which we denote by $\mbox{Sym}(X)$.  This construction satisfies the adjointness relationship 
$$ Hom_{\mbox{Alg}_S}(\mbox{Sym}(X), A) \cong Hom_{\mbox{Mod}_S}(X, A)
$$
 Sym is a triple on the category of  spectra, and commutative  $S$-algebras are exactly the algebras over the triple
Sym.   For any commutative  $S$-algebra $A$ , we may as in Definition \ref{simpresolution} construct its simplicial resolution $\mbox{Sym}_{\cdot}(A)$, relative to the triple Sym.   Moreover, Proposition \ref{A} shows that the natural map $\mbox{Sym}_{\cdot}(A) \rightarrow A$ is a weak equivalence.  In \cite{shipley}, it is shown that $\mbox{Sym}_{\cdot}(A)$ is a  cofibrant object in the model structure on the category $\mbox{Alg}_S$.  

\noindent Given two commutative  $S$-algebras $A$ and $B$, the set $Hom_{\mbox{Alg}_S}(A,B)$ is  equipped with the structure of a space, as a subspace of the zeroth space of the spectrum $Hom_{\mbox{Mod}_S}(A,B)$.  This construction is not homotopy invariant, but  the result of replacing $A$ by any weakly equivalent cofibrant object does yield a homotopy invariant notion.  Of course, we may take the cofibrant replacement to be $\mbox{Sym}_{\cdot}(A)$.  We now need a result about the space of $S$-algebra maps to Eilenberg-MacLane spectra.  

\begin{Proposition} Let $A$ be a (-1)-connected commutative  $S$-algebra.  Let B denote any commutative ring, and $\Bbb{H}(B)$ the corresponding Eilenberg-MacLane spectrum.  The canonical map from the space $Hom_{\mbox{ \em Alg}_S}(\mbox{Sym}_{\cdot}(A),\Bbb{H}(B)) \rightarrow 
Hom(\pi _0 (A), B)$ is an equivalence, i.e. $Hom_{\mbox{ \em Alg}_S}(\mbox{Sym}_{\cdot}(A),\Bbb{H}(B))$ is a space whose components are in bijective correspondence with the ring homomorphisms from $\pi _0 (A)$ to $B$, and so that each component is contractible.  
\end{Proposition}
\begin{Proof} The adjunction $ Hom_{\mbox{Alg}_S}(\mbox{Sym}(X), A) \cong Hom_{\mbox{Mod}_S}(X, A)
$ shows that the space $$Hom_{\mbox{Alg}_S}(\mbox{Sym}_{\cdot}(A),\Bbb{H}(B))$$  is equivalent to the total space of the cosimplicial space which in level $k$ is the space of spectrum maps from $\mbox{Sym}^k(A)$ to $\Bbb{H}(B)$.  Since the  $k+1$-skeleton   $Sk^{(k+1)}\mbox{Sym}_{\cdot}(A) $ is obtained from the $k$-skeleton by attaching cells in dimensions $k+1$ and higher,  it follows that the restriction map 
$$ Hom_{\mbox{Alg}_S}(\mbox{Sym}_{\cdot}(A),\Bbb{H}(B)) \rightarrow Hom_{\mbox{Alg}_S}(Sk^1 \mbox{Sym}_{\cdot}(A),\Bbb{H}(B))\stackrel{defn}{=} T^1
$$
is an equivalence.  Now, $T^1$ is clearly the homotopy equalizer of the two maps $$\varphi, \psi : Hom_{\mbox{Mod}_S}(A, \Bbb{H}(B)) \rightarrow  Hom_{\mbox{Mod}_S}(\mbox{Sym}(A), \Bbb{H}(B))$$ defined by 
$$ \left\{
\begin{array}{l}
\varphi (f) = f \compcirc \alpha _A \mbox{ where $\alpha _A: \mbox{Sym}(A) \rightarrow A$ is the structure map for $A$} \\
\psi (f) = \alpha _B \compcirc \mbox{Sym}(f) \mbox{ where $\alpha _{\Bbb{H}(B)}: \mbox{Sym}(\Bbb{H}(B))\rightarrow \Bbb{H}(B)$ is the structure map for $\Bbb{H}(B)$.} 
\end{array}
\right . 
$$
The space $Hom_{\mbox{Mod}_S}(A, \Bbb{H}(B))$ is clearly equivalent to the discrete space $Hom_{Ab}(\pi _0 (A),  B)$.  Similarly,   $Hom_{\mbox{Mod}_S}(\mbox{Sym}(A), \Bbb{H}(B))$ is  equivalent to the discrete space $$Hom_{Ab}(\pi _0 (\mbox{Sym}(A)), B) = Hom_{Ab}(\mbox{Sym}(\pi _0 (A), B))$$
The homotopy equalizer is now equivalent to the equalizer of the pair of set maps 
$$\Phi, \Psi:   Hom_{Ab}(\pi _0 (A),  B) \rightarrow Hom_{Ab}(\mbox{Sym}(\pi _0 (A), B))
$$
defined by 
$$
\left\{ 
\begin{array}{l}
\Phi (f) = f \compcirc \alpha _{\pi _0A} \mbox{ where $\alpha _{\pi _0 A}: \mbox{Sym}(\pi _0 A) \rightarrow \pi _0 A$ is the structure map for the ring $\pi _0 A$.} \\
\Psi (f) = \alpha _B \compcirc \mbox{Sym}(f) \mbox{ where $\alpha _B : \mbox{Sym}(B) \rightarrow B$ is the structure map for the ring $B$}. 
\end{array}
\right . 
$$
This set is the set of all abelian group homomorphisms $f : \pi _0 A \rightarrow B$ making the diagram 
$$
\begin{diagram} \node{\mbox{Sym}(\pi _0 A)} \arrow{s,l}{\alpha _{\pi _0 A}} \arrow{e,t}{\mbox{Sym}(f)} \node{\mbox{Sym}(B)} \arrow{s,r}{\alpha _B}\\
\node{\pi _0 A} \arrow{e,t}{f} \node{B}
\end{diagram}
$$
commute.  This is clearly the set of ring homomorphisms from $\pi _0 A$ to $B$. 
\end{Proof}
\begin{Corollary} There is a canonical homotopy equivalence class of homomorphisms of commutative  $S$-algebras $\pi _A : A \rightarrow \Bbb{H}(\pi _0A)$ which induces the identity on $\pi _0$. 
\end{Corollary}
We also have the following results on $A$-module spectra.  They follow immediately from  Proposition 3.9 of \cite{DGI}.  

\begin{Proposition} \label{UC}  Let $A$ be a (-1)-connected commutative  $S$-algebra, and let $M$ be a cofibrant $(k-1)$-connected $A$-module spectrum.  Let $N$ be any module over the ring $\pi _0 A$, and $\Bbb{H}(N,k)$ the $k$-dimensional Eilenberg-MacLane spectrum for $N$,  regarded as an $A$-module spectrum via the homomorphism $A \rightarrow \pi _0 A$.  The the zeroth space of the function spectrum $Hom_A( M, \Bbb{H}(N, k))$ is homotopy equivalent to the set $Hom_{\pi _0 A}(\pi _k M, N)$ via the natural map which assigns to any map its induced map on $\pi _k$.   
\end{Proposition}

\begin{Corollary} \label{classify} Let $A$ be a (-1)-connected commutative $S$-algebra.  Let $M$ be a  module spectrum over $A$, and suppose that the underlying spectrum of $M$ is an Eilenberg-MacLane spectrum, say of dimension $n$.  $\pi _n M$ is  a module over the ring $\pi _0 A$. We let $\Pi _M$ denote the $A$-module spectrum obtained from $\pi _n M$ by pullback along the ring homomorphism $\pi _A$.  Then, $M$ is equivalent as an $A$-module spectrum to $\Pi _M$.  
\end{Corollary}

\noindent These results allow us to construct a Postnikov tower in the category of connective modules over $A$. 
\begin{Corollary}Let $A$ be a (-1) connected commutative  $S$-algebra, and let $M$ be a $k$-connected $A$-module spectrum, for some integer $k$.  Then there is a tower of fibrations of $A$-module spectra
$$  \cdots M[s] \rightarrow M[s] \rightarrow \cdots \rightarrow M[k+1] \rightarrow M[k]
$$
together with $A$-module spectrum maps $M \stackrel{f_s}{\rightarrow} M[s]$ making all the triangles 
$$
\begin{diagram}
\node{M} \arrow{e,t}{f_s} \arrow{se,b}{f_{s-1}} \node{M[s]} \arrow{s} \\
\node{} \node{M[s-1]}
\end{diagram}
$$
commute, with the following properties. 
\begin{itemize}
\item{$f_s$ induces an isomorphism on $\pi _i$ for $i \leq s$. }
\item{$\pi _i M[s] = 0 $ for $i > s$. }
\item{The homotopy fiber of the map $M[s] \rightarrow M[s-1]$ is an Eilenberg-MacLane spectrum in dimension $s$. }
\item{ Using the maps $\{ f_s \}_s$, $\holim{s}{M[s]}$ is naturally equivalent to $M$}
\end{itemize} 
Moreover, the tower is natural in the homotopy category. 
\end{Corollary}
\begin{Proof} Immediate from Proposition \ref{classify} above.  
\end{Proof}

\noindent We will use the Postnikov tower to construct the spectral sequence.  A preliminary observation is the following. 

\begin{Proposition} Let 
$$ \cdots \rightarrow M_s \rightarrow \cdots \rightarrow M_1 \rightarrow M_0
$$
be any diagram of module spectra over a commutative  $S$-algebra $A$, and let $P$ denote any cofibrant $A$-module spectrum.  Suppose further that the connectivity of the $M_s$'s goes to infinity with $s$, and that $P$ is connective, i.e. is $k$-connected for some $k$.  Then the natural map 
$$ P \smashover{A}( \holim{s}{M_s}) \longrightarrow  \holim{s}{(P \smashover{A}M_s)}
$$
is an equivalence of spectra.  
\end{Proposition}
\begin{Proof} Follows directly from the fact that $P\smashover{A}$ preserves connectivity when $P$ is cofibrant. 
\end{Proof}
\begin{Corollary} \label{converge} Suppose that $A$ is a (-1)-connected commutative  $S$-algebra, and that we have a homomorphism $f: A \rightarrow B$ of commutative  $S$-algebras.  Suppose further that $B$ is (-1)-connected and cofibrant.  For any diagram 
$$ \cdots \rightarrow M_s \rightarrow \cdots \rightarrow M_1 \rightarrow M_0
$$
of connective $A$-module spectra, for which the connectivity of $M_s$ goes to infinity with $s$, the natural map 
$$ (\holim{s}{M_s})^{\wedge}_B \rightarrow \holim{s}{(M_s)^{\wedge}_B}
$$
is a weak equivalence of spectra. In particular, if $\{ M[s] \} _s$ is the Postnikov tower for an $A$-module spectrum $M$, then $M^{\wedge}_B \cong \holim{s}{M[s]^{\wedge}_B}$.  Further, the homotopy fiber of the map $M[s]^{\wedge}_B \rightarrow M[s-1]^{\wedge}_B$ is equivalent to $(F_s)^{\wedge}_B$, where $F_s$ is the homotopy fiber of the map $M[s] \rightarrow M[s-1]$. 
\end{Corollary}

\noindent The spectral sequence in question is simply the spectral sequence on homotopy groups attached to an inverse system of spectra.  This is a standard result which is discussed in \cite{Bousfield}, Chapter IX. 

\begin{Proposition} \label{generic}  Let $ \cdots X_s \rightarrow X_{s-1} \rightarrow \cdots \rightarrow X_1 \rightarrow X_0$ be an inverse system of connective spectra,  and let $F_s$ denote the homotopy fiber of the map $X_s \rightarrow X_{s-1}$.  Suppose further that the connectivity of the  spaces goes to infinity with $s$.  Then  there exists  a left half plane spectral sequence  with $E_1^{pq}$-term $\pi _{p+q}F_{-p}$, converging to $\pi _{p+q}(\holim{s}{X_s})$. 
\end{Proposition}

\noindent {\bf Proof of Theorem \ref{algebraictogeometric}:} We apply Proposition \ref{generic} and Corollary \ref{generic} to the inverse system 

$$  \cdots \rightarrow M[s]^{\wedge}_B \rightarrow M[s-1]^{\wedge}_B \rightarrow \cdots \rightarrow M[k+1]^{\wedge}_B \rightarrow M[k]^{\wedge}_B
$$
where $k$ is the connectivity of $M$.   The homotopy fiber of the map $M[s]^{\wedge}_B \rightarrow M[s-1]^{\wedge}_B$ is $(F_s)^{\wedge}_B$, where as above $F_s$ denotes the homotopy fiber of the map $M[s] \rightarrow M[s-1]$.  Since $\{ M[s] \}_s$ is the Postnikov tower for $M$, $F_s$ is an $s$-dimensional Eilenberg-MacLane spectrum with $\pi _s (F_s) \cong \pi _s (M)$.  Consequently, by Corollary \ref{classify}, it is equivalent to a module over the ring $\pi _0(A)$, and by Theorem \ref{ringinvar}, it follows that the derived completion of $F_s$ at  the $S$-algebra  homomorphism $A \rightarrow B$ is equivalent to its derived completion at the ring homomorphism $\pi _0 A \rightarrow \pi _0 B$.  The theorem now follows. \qed

\noindent We can now use this spectral sequence to obtain two useful corollaries. We first have an extension of Proposition \ref{properties}, item  6. 

\begin{Corollary}  Let $B \rightarrow C$ be a homomorphism of (-1)-connected commutative $S$-algebras, and suppose $\pi _0B \rightarrow \pi _0 C $ is surjective.  Let $M$ be a $B$-module spectrum, and suppose that for each $i$, the $\pi _0B$-module structure on $\pi _1 M$ extends to  a $\pi _0 C$-action.  Then the natural map $ \eta : M \rightarrow M^{\wedge}_C$ is an equivalence.  
\end{Corollary}
\begin{Proof} Follows directly from Proposition \ref{properties}, item 6, and Theorem \ref{algebraictogeometric}. 
\end{Proof}

\noindent We also have the following corollary, which will be useful in the $K$-theoretic applications.  

\begin{Corollary} Let 
$$ A \longrightarrow B \longrightarrow C
$$
be a diagram of (-1)-connected commutative $S$-algebras, and suppose that $\pi _0 A \rightarrow \pi _0 B$ and $\pi _0 B \rightarrow \pi _0 C$  are surjective.  For any $B$-module $M$, let $\rho _A(M)$ denote $M$ regarded as an $A$-module by restriction of scalars.  Then the natural map $(\rho _A(M))^{\wedge}_C \rightarrow M^{\wedge}_C$ is a weak equivalence.  
\end{Corollary}
\begin{Proof} It follows from Theorem \ref{algebraictogeometric} that it suffices to consider the case where $A,B$, and $C$ are all obtained by applying the Eilenberg-MacLane construction to diagram of actual commutative rings, and that $M$ is obtained as an actual module over the ring $B$.  Exactly as in the proof of Theorem \ref{ringinvar}, it suffices to prove that the natural map 
$$  C \smashover{A} M \longrightarrow (C \smashover{A} M) ^{\wedge}_C
$$ 
is an equivalence.  Here the completion is taking place along the homomorphism $B \rightarrow C$, and the $B$-module structure on $C \smashover{A} M$ is obtained from the original $B$-action on $M$, not via restriction of scalars from the extended left $C$-action.  That is, we have $b \cdot (c \wedge  m) = c \wedge b \cdot m$. Now, the homotopy groups of $C \smashover{A} M$ are the groups $Tor^A_* (M,C)$, via the (collapsing) K\"{u}nneth  spectral sequence.  We must only show that the $B$ action on these $Tor$-groups obtained from the given $B$-action on $M$ extends over $C$.  But from the surjectivity of the homomorphism $A \rightarrow B$, the $B$-action can be computed as the action of $A /J$, where $J$ is the kernel of the ring homomorphims $A \rightarrow B$, and hence is determined by the $A$-action.  It is clear from the definition that the $A$-action extends over $C$, and hence that it vanishes on the $I$, the kernel of the ring homomorphism $A \rightarrow C$.  But, $A/I \cong B/I^{\prime}$, where $I^{\prime}$ is the kernel of $B \rightarrow C$, which gives the result.   \end{Proof}

  \section{Examples}
  
  \noindent In this section, we will discuss a number of examples of this construction.   
  
  \begin{Example} {\em We consider the case where the $S$-algebras $A$ and $B$  are the Eilenberg-MacLane spectra for the rings $\Bbb{Z}$ and $\Bbb{F}_p$, and where the homomorphism $f:A \rightarrow B$ is induced by reduction mod $p$.  For a finitely generated $A$-module $M$,  Theorem \ref{ringcase} tells us that the derived completion of $M$ along $f$ is simply the Eilenberg-MacLane spectrum for the completed module $M^{\wedge}_p$.  }\end{Example} 
  \begin{Example}{ \em Let $A= S^0$, the sphere spectrum, and as in the previous example let $B$ be the Eilenberg-MacLane spectrum for the ring $\Bbb{F}_p$.  Let $f$ be the composition of the natural map the the Eilenberg-MacLane spectrum for $\Bbb{Z}$ with mod-$p$ reduction.  Any spectrum is in a natural way a module over the commutative  $S$-algebra $S^0$, and its completion along this map is the usual $p$-adic completion.  The tower of fibrations mentioned in Remark \ref{filtered} is in this case identical to the tower in the Adams spectral sequence.  } 
  \end{Example}
  \begin{Example}{\em Let $A$ be the $S$-algebra $ku$, i.e. connective complex $K$-theory, and let $B$ be the Eilenberg-MacLane spectrum for the ring $\Bbb{Z}$.  We view $ku$ as the spectrum associated to the topological symmetric monoidal category of finite dimensional complex vector spaces, and we define $f$ to be the homomorphism of  $S$-algebras induced by the functor which sends every complex vector space to its dimension.  Since $\pi _0 f$ is an isomorphism, Theorem \ref{algebraictogeometric} shows that in this case, $ku^{\wedge}_B$ is equivalent to $ku$ itself.  The tower of fibrations mentioned in Remark \ref{filtered} is in this case the {\em accelerated Postnikov tower} for $ku$, with the $i$-th spectrum in the tower of fibrations equivalent to the Postnikov cover $ku[2i, \ldots , + \infty )$}
  \end{Example}
  \begin{Example} {\em  Let $A$ denote the Eilenberg-MacLane spectrum for the representation ring $R[\Bbb{Z}_p]$, where $\Bbb{Z}_p$ denotes the group of $p$-adic integers.  For a profinite group $G$, the representation ring $R[G]$ is defined to be the direct limit of $R[G/N]$, where $N$ ranges over all the finite index normal subgroups of the group $G$.  As a ring, one readily sees that $R[\Bbb{Z}_p] \cong \Bbb{Z}[ \Bbb{Z}/{p^{\infty}}\Bbb{Z}]$, where $\Bbb{Z} [ - ] $ denotes integral group ring, and where $ \Bbb{Z}/{p^{\infty}}\Bbb{Z}$ denotes the union $\bigcup _n \Bbb{Z}/p^n \Bbb{Z}$.  We let $C = \Bbb{F}_p$, and define $f:A \rightarrow C$ to be the composite $R[\Bbb{Z}_p] \rightarrow R[\{ e \} ] = \Bbb{Z} \rightarrow \Bbb{F}_p$.  We will let the module $M$ be the ring $A$ itself.  We let $S^1_{\cdot}$ denote the simplicial group of singular chains on the circle group,  with multiplication and inverse induced by those on the topological group. There is of course a natural homomorphism of simplicial groups $ \Bbb{Z}/{p^{\infty}}\Bbb{Z} \rightarrow S^1_{\cdot}$, induced by the identification of $ \Bbb{Z}/{p^{\infty}}\Bbb{Z}$ with the group of $p$-power torsion points on $S^1_{\cdot}$, and where $ \Bbb{Z}/{p^{\infty}}\Bbb{Z}$ is regarded as a discrete simplicial group. 
  We define $B$ to be the simplicial group ring $\Bbb{Z}[S^1_{\cdot} ]$, which is a simplicial  algebra  $ \Bbb{Z}[\Bbb{Z}/{p^{\infty}}\Bbb{Z}]$.  Such a simplicial algebra may be regarded as a commutative $S$-algebra, in fact a commutative algebra spectrum  over the $S$-algebra $\Bbb{H}(A)$.  We let $\rho (B)$ denote $B$ regarded as a $A$-module spectrum via restriction of scalars along the inclusion $\Bbb{H}(A) \rightarrow B$.  We then have a composite
  
\begin{equation} \label{comparison} \Bbb{H}(A) _C^{\wedge} \rightarrow \rho (B)^{\wedge}_C \rightarrow B^{\wedge}_C
\end{equation}
  Note that the right  hand spectrum is obtained by completion over the ground $S$-algebra $B$, and the two leftmost spectra are completions over the ground $S$-algebra $\Bbb{H}(A)$. The right hand map is induced by the natural morphism of cosimplicial spectra ${\cal T}^{\cdot}_{\Bbb{H}(A)}(\rho (B);C) \rightarrow {\cal T}^{\cdot}_{B} (B;C)$, so the composite is induced by the natural morphism ${\cal T}^{\cdot}_{\Bbb{H}(A)} (A;C) \rightarrow {\cal T}^{\cdot}_B(B;C) $.  One can now check from the definitions that 
  $$  \pi _s ({\cal T}^{k}_A (A;C)) \cong \mbox{MultiTor}^A_{s}(\underbrace{\Bbb{F}_p, \Bbb{F}_p, \ldots , \Bbb{F}_p}_{k+1 \mbox{ factors }}) \cong H_s( \underbrace{B\Bbb{Z}/p^{\infty}\Bbb{Z} \times \cdots \times B\Bbb{Z}/p^{\infty}\Bbb{Z}}_{k \mbox{ factors}}, \Bbb{F}_p)
  $$
  and
  $$  \pi _s ({\cal T}^{k}_B (B;C)) \cong H_s( \underbrace{BS^1 \times \cdots \times BS^1}_{k \mbox{ factors}}, \Bbb{F}_p)
$$  It is easy to check that the induced homomorphism $ \pi _s ({\cal T}^{k}_A (A;C)) \rightarrow  \pi _s ({\cal T}^{k}_B (B;C))$ is (under the isomorphisms above) identified with $H_s( \underbrace{Bi \times \cdots \times Bi}_{k \mbox{ factors} },  \Bbb{F}_p)$, where $i \colon \Bbb{Z}/p^{\infty}\Bbb{Z} \rightarrow S^1$ is the inclusion. It is a well known fact that this map is an isomorphism, which shows that the map in \ref{comparison} above is a an equivalence of spectra.  In order to identify  $ \Bbb{H}(A) _C^{\wedge}$,  then, it will suffice to identify $B^{\wedge}_C$.  To do this, we consider the homomorphism of commutative $S$-algebras $g: S^0 \rightarrow \Bbb{H}(Z) \rightarrow \Bbb{Z}[S^1]$.  As before, let $\rho (\Bbb{Z}[S^1])$ denote $\Bbb{Z}[S^1]$ regarded as an $S^0$-module via restriction of scalars along $g$.  Then we obtain a natural map $\rho ((\Bbb{Z}[S^1])^{\wedge}_{\Bbb{H}(\Bbb{F}_p)} \rightarrow \Bbb{Z}[S^1]^{\wedge}_{\Bbb{H}(\Bbb{F}_p)}$. Again, the left hand completion is over the commutative $S$-algebra $S^0$, and the right one is over the $S$-algebra $\Bbb{Z}[S^1]$.  There is a homomorphism of algebraic to geometric spectral sequences (\ref{algebraictogeometric}), which is an isomorphism since $\pi _0S^0 \rightarrow \pi _0 \Bbb{Z}[S^1]$ is.  Consequently, the desired completion is the $p$-adic completion of the spectrum $\Bbb{Z}[S^1]$.  The resulting homotopy groups are the groups $H_i (S^1, \Bbb{Z}_p)$, so we finally obtain the formula 
$$\begin{array}{l}{\pi  _i (A^{\wedge}_{\Bbb{H}(\Bbb{F}_p)} )\cong \Bbb{Z}_p \mbox{ for }i=0,1} \\
{\pi  _i (A^{\wedge}_{\Bbb{H}(\Bbb{F}_p)}) = 0 \mbox{ \hspace{.1cm} otherwise}}
\end{array}$$

\noindent This result has been extended to a statement about finitely generated nilpotent groups by T. Lawson in \cite{lawson}. A construction referred to as {\em deformation $K$-theory} can be made which incorporates the topology on spaces of representations, and which applies to any discrete group.  In the case of $\Bbb{Z}$, it produces a spectrum with homotopy groups given by $\Bbb{Z}[S^1_.] \otimes \pi _* ku$.  What Lawson proves is that for a finitely generated nilpotent group $\Gamma$, the derived completion of the representation ring of the pro-$p$ completion of $\Gamma$ at the homomorphism given by mod $p$ augmentation is equivalent to the $p$-adic completion of the deformation $k$-theory of $\Gamma$.  This gives a link between this completion process on representation rings of pro-$p$ quotients of $\Gamma$ with the geometry of representation varieties of $\Gamma $. }

  \end{Example}

 \end{document}